\documentclass[12pt]{article}

\usepackage{amsmath}
\usepackage{amsfonts}
\usepackage{amssymb}
\usepackage{latexsym}
\usepackage{graphicx}

\usepackage{amsmath, amssymb,latexsym}
\usepackage{color}
\usepackage{fullpage}

\newtheorem{theorem}{Theorem}[section]
\numberwithin{equation}{section}

\DeclareMathOperator{\diag}{diag}

\DeclareMathOperator{\SL}{SL}
\DeclareMathOperator{\GL}{GL}
\DeclareMathOperator{\SO}{SO}
\DeclareMathOperator{\Oo}{O}
\DeclareMathOperator{\CO}{CSO}
\DeclareMathOperator{\U}{U}
\DeclareMathOperator{\Curl}{Curl}
\newcommand{\A}{Z} 
\newcommand{\id} {I}
\newcommand{\logp} {\log}
\DeclareMathOperator{\Log}{Log}

\renewcommand{\skew}{\mathop{\rm skew_{\!_*}\!}}
\DeclareMathOperator{\rskew}{skew}
\DeclareMathOperator{\sym}{sym_{_*}\!}
\DeclareMathOperator{\rsym}{sym}

\DeclareMathOperator{\Tr}{tr}

\DeclareMathOperator{\dev}{dev}
\DeclareMathOperator{\sL}{\mathfrak{sl}}
\DeclareMathOperator{\so}{\mathfrak{so}}

\DeclareMathOperator{\Det}{det}
\renewcommand{\det}[1]{ {\Det({#1})} }
\newcommand{\tr}[1]{ {\Tr \left({#1}\right)} }
\newcommand{\ds}{\,{\rm{ds}} }

\newcommand{\norm}[1]{\|#1\|}
\newcommand{\R}{\mathbb{R}}

\newcommand{\C}{\mathbb{C}}
\newcommand{\N}{\mathbb{N}}

\renewcommand{\Re}{\mathop{\mathfrak{Re}} }
\renewcommand{\Im}{\mathop{\mathfrak{Im} }}

\DeclareMathOperator{\dist}{dist}

\DeclareMathOperator{\arcosh}{arcosh}

\newcommand{\abs}[1]{|#1|}
\newcommand{\Mprod}[2]{ {\langle #1 ,#2\rangle} }

\newcommand{\qeda}{\tag*{$\blacksquare$} }

\title{
A logarithmic minimization property of the unitary polar factor in the spectral norm and the Frobenius matrix norm
}

\author{Patrizio Neff%
\thanks{Corresponding author, 
Head of Lehrstuhl f\"ur Nichtlineare Analysis und Modellierung, Fakult\"at f\"ur  Mathematik, Universit\"at Duisburg-Essen, Campus Essen, Thea-Leymann Str. 9, 45127 Essen, Germany, email: patrizio.neff@uni-due.de, Tel.: +49-201-183-4243}\; \and 
Yuji Nakatsukasa\thanks{
Department of Mathematical Informatics, University of Tokyo, Tokyo 113-8656, Japan, email: nakatsukasa@mist.i.u-tokyo.ac.jp}\; 
\and
Andreas Fischle\thanks{
Fakult\"at f\"ur  Mathematik, Universit\"at Duisburg-Essen, Campus Essen, Thea-Leymann Str. 9, 45127 Essen, Germany, email: andreas.fischle@uni-due.de}%
}

\begin{document}
\maketitle

\begin{abstract}
  The unitary polar factor $Q=U_p$ in the polar decomposition of $\A=U_p \, H$ is the minimizer over unitary matrices $Q$ for both $\|\Log(Q^* \A)\|^2$ and its Hermitian part $\|\sym(\Log(Q^* \A))\|^2$ over both $\mathbb{R}$ and $\mathbb{C}$ for any given invertible matrix $\A\in\C^{n\times n}$ and any matrix logarithm $\Log$, not necessarily the principal logarithm $\log$. We prove this  for the spectral matrix norm for any $n$ and for the Frobenius matrix norm in two and three dimensions. The result shows that the unitary polar factor is the nearest orthogonal matrix to $Z$ not only in the normwise sense, but also in a geodesic distance. The derivation is based on Bhatia's generalization of Bernstein's trace inequality for the matrix exponential and a new sum of squared logarithms inequality. Our result generalizes the fact for scalars that for any complex logarithm and for all $z\in\C\setminus\{0\}$
\begin{align*}
 \min_{\vartheta\in(-\pi,\pi]}\abs{\Log_{\C}(e^{-i\vartheta} z)}^2&=\abs{\log\abs{z}}^2\, ,\quad \min_{\vartheta\in(-\pi,\pi]}\abs{\Re\Log_{\C}(e^{-i\vartheta} z)}^2=\abs{\log\abs{z}}^2\, .
\end{align*}
\end{abstract} 



\section{Introduction}
Every matrix $\A\in\mathbb{C}^{m\times n}$ admits a polar decomposition 
\begin{align*}\A=U_p\, H\, ,\end{align*}
where the \emph{unitary polar factor} $U_p$ has orthonormal columns and $H$ is Hermitian positive semidefinite~\cite{Autonne1902}, \cite[Ch.~8]{Higham2008}. The decomposition is unique if $\A$ has full column rank. In the following we assume that $\A$ is an invertible matrix, in which case $H$ is positive definite. The polar decomposition is the matrix analog of the polar form of a complex number 
\begin{align*}
   z= e^{i\, \arg(z)}\cdot r, \quad r=\abs{z}\ge 0 ,\quad -\pi< \arg(z)\le \pi\, .
\end{align*}
The polar decomposition has a wide variety of applications:  the solution to the Euclidean orthogonal Procrustes problem $\min_{Q\in \U(n)}\|\A-BQ\|_F^2$ is given by the unitary polar factor of $B^*\A$~\cite[Ch.~12]{Golubbookori}, and the polar decomposition can be used as the crucial tool for computing the eigenvalue decomposition of symmetric matrices and the singular value decomposition (SVD)~\cite{nakahig12}. Practical methods for computing the polar decomposition are the scaled Newton iteration \cite{ByXu08} and the QR-based dynamically weighted Halley iteration~\cite{nakatsukasa:2700}, and their backward stability is established in~\cite{nahi11}. 

The unitary polar factor $U_p$ has the important property~\cite[Thm.~IX.7.2]{Bhatia97}, \cite{Fan55}, \cite[p.~197]{Higham2008}, \cite[p.454]{Horn85} that it is the nearest unitary matrix to $\A\in\C^{n\times n}$, that is, 
\begin{align}
\label{unitary_polar_optimal}
   \min_{Q\in\U(n)}\norm{\A-Q}^2= \min_{Q\in\U(n)}\| Q^* \A-\id\|^2=\| U_p^* \A-\id\|^2=\norm{\sqrt{\A^*\A}-\id}^2\, ,
 \end{align}
 where $\| \cdot\|$ denotes \emph{any unitarily invariant norm}.
For the Frobenius matrix norm in the three-dimensional case the proof of this optimality theorem was first given by Grioli \cite{Grioli40} in an article on the theory of elasticity, an annotated translation of which can be found in \cite{OnGrioli}.
The optimality for the Frobenius norm implies for real $\A\in\R^{n\times n}$ and the orthogonal polar factor \cite{Guidugli80}
\begin{align}
\label{optimal_real_trace}
    \forall\; Q\in\Oo(n):\quad  \tr{Q^T \A}= \Mprod{Q}{\A}\le \Mprod{U_p}{\A}=\tr{U_p^T \A}\, .
 \end{align}
 In the complex case we similarly have \cite[Thm.~7.4.9, p.432]{Horn85}
  \begin{align}
    \forall\; Q\in\U(n):\quad  \Re \tr{Q^* \A} \le 
\Re \tr{U_p^* \A}\, .
 \end{align} 
 For invertible $\A\in\GL^+(n,\R)$ and the Frobenius matrix norm $\norm{\cdot}_F$ it can be shown that~\cite{Pietraszkiewicz05,Guidugli80} 
 \begin{align}
   \label{polar_classic}
    \min_{Q\in \Oo(n)} \left(\mu\, \| \sym (Q^T \A-\id)\|_F^2+\mu_c\, \|\skew (Q^T \A-\id)\|_F^2\right)=\mu\, \| U_p^T \A-\id\|_F^2 \, ,
 \end{align} 
 for $\mu_c\ge \mu>0$.  Here, $\sym (X)=\frac{1}{2}(X^*+X)$ is the Hermitian part and $\skew (X)=\frac{1}{2}(X-X^*)$ is the skew-Hermitian part of $X$. The family \eqref{polar_classic} appears as the strain energy expression in geometrically exact Cosserat extended continuum 
 models \cite{Klawonn_Neff_Rheinbach_Vanis09,Neff_plate04_cmt,Neff_micromorphic_rse_05,Neff_Cosserat_plasticity05,Neff_Muench_simple_shear09}. 

Surprisingly, the optimality \eqref{polar_classic} of the orthogonal polar factor ceases to be true for $0\le \mu_c<\mu$. Indeed, for $\mu_c=0$ 
 there exist $\A\in\R^{3\times 3}$ such that $U_p^T$ is not even a minimizer in the special orthogonal group~\cite{Neff_Biot07}:
\begin{align}
 \label{counter_example}
 \min_{Q\in\Oo(3)} \| \rsym (Q^T \A-\id)\|_F^2 \:\leq \min_{Q\in\SO(3)} \| \rsym (Q^T \A-\id)\|_F^2 \,<\, \|U_p^T \A-\id\|_F^2\, .
 \end{align} 
Note that we drop the subscript $_\ast$ in $\rm sym_\ast$ and $\rm skew_\ast$ when the matrix is real. By compactness of $\Oo(n)$ and $\SO(n)$ and the continuity of $Q\mapsto \norm{\rsym Q^T\A-\id}_F^2$ it is clear that the minima in \eqref{counter_example} exist.
Here, the polar factor $U_p$ of $\A\in\GL^+(3,\R)$ is always a critical point, but is not necessarily (even locally) minimal. In contrast to $\norm{X}_F^2$ the term $\norm{\rsym X}_F^2$ is not invariant w.r.t. left-action of $\SO(3)$ on $X$, which does explain the appearance of nonclassical solutions in \eqref{counter_example} since now
 \begin{align}
 \label{symmetry_break}
   \norm{\rsym Q^T\A-\id}_F^2=\norm{\rsym Q^T\A}_F^2-2\, \tr{Q^T \A}+3\, 
 \end{align}
and optimality does not reduce to optimality of the trace term \eqref{optimal_real_trace}. The reason that there is no nonclassical solution in \eqref{polar_classic} is that for $\mu_c\ge \mu>0$ we have
\begin{align}
 \min_{Q\in \Oo(n)} &\mu\, \| \rsym (Q^T \A-\id)\|_F^2+\mu_c\, \|\rskew (Q^T \A-\id)\|_F^2  \\
 \ge &\min_{Q\in \Oo(n)} \mu\, \| \rsym (Q^T \A-\id)\|_F^2+\mu\, \|\rskew (Q^T \A-\id)\|_F^2
 =  \min_{Q\in \Oo(n)} \mu\, \| Q^T \A-\id\|_F^2\, . \notag
 \end{align}

\subsection{The matrix logarithm minimization problem and results}\label{subsec:minproblem}
Formally, we obtain our minimization  problems 
 \[\min_{Q\in\SO(n)}  \|  \Log(Q^T \A)\|_F^2 \quad\text{ and } \quad\min_{Q\in\SO(n)}  \| \sym \Log(Q^T \A)\|_F^2\]
by replacing the matrix $Q^T \A-\id$ by the matrix $\Log(Q^T \A)$
in \eqref{polar_classic}. Then, introducing the weights $\mu,\mu_c\ge 0$ we embed the problem in a more general family of minimization problems at a given $\A\in \GL^+(n,\R)$
  \begin{align}
    \min_{Q\in\SO(n)} \mu\, \| \sym \Log(Q^T \A)\|_F^2+\mu_c\, \|\skew \Log(Q^T \A)\|_F^2\, ,\quad \mu>0,\, \mu_c\ge 0\, .
    \label{general_problem}
 \end{align} 
For the solution of \eqref{general_problem} we consider separately the minimization of 
  \begin{align}
   \label{problem2}
    \min_{Q\in \U(n)} \| \Log(Q^* \A)\|^2\, ,\quad \min_{Q\in \U(n)} \| \sym\Log(Q^* \A)\|^2\, ,
    \end{align}  
 on the group of  \emph{unitary matrices} $Q\in\U(n)$ and with respect to \emph{any matrix logarithm} $\Log$.
 
We show that the unitary polar factor $U_p$ is a minimizer of both terms in \eqref{problem2} for both the Frobenius norm (dimension $n=2,3$) and the spectral matrix norm for arbitrary $n\in\N$, and the minimum is attained when the principal logarithm is taken.
 
 Finally, we show that the minimizer of the real problem \eqref{general_problem} for all $\mu>0,\; \mu_c\ge 0$ is also given by the polar factor $U_p$. Note that $\sym (Q^T\A-\id)$ is the leading order approximation to the Hermitian part of the logarithm $\sym\Log Q^T\A$ in the neighborhood of the identity $\id$, and recall the non-optimality of the polar factor in \eqref{counter_example} for $\mu_c=0$. The optimality of the polar factor in \eqref{problem2} is therefore rather unexpected. Our result implies also that different members of the family of Riemannian metrics $g_X$ on the tangent space \eqref{Riemannian_metric} lead to the same Riemannian distance  to the compact subgroup $\SO(3)$, see \cite{Neff_Osterbrink_hencky13}. 
  
 Since we prove that the unitary polar factor $U_p$ is the unique minimizer for \eqref{general_problem} and the first term in \eqref{problem2}
for the Frobenius matrix norm and $n\leq 3$, it follows that these new optimality properties of $U_p$ provide another  characterization of the polar decomposition.

In our optimality proof we do not use differential calculus on the nonlinear manifold $\SO(n)$ for the real case because the derivative of the matrix logarithm is 
analytically not really tractable. However, if we assume  a priori that the minimizer $Q_\sharp\in\SO(n)$ can be found in the set $\{Q\in\SO(n)\,|\, \norm{Q^T \A-\id}_F\le q<1\,\}$, we can use the power series expansion of the principal logarithm and differential calculus to show that the polar factor is indeed the unique minimizer (we hope to report this elsewhere). 

Instead, motivated by insight gained in the more readily accessible complex case, we first consider the Hermitian minimization problem, which has the advantage of allowing us to work with the positive definite Hermitian matrix $\exp[\sym\Log Q^* \A]$. A subtlety that we encounter several times is the possible non-uniqueness of the matrix logarithm $\Log$. 
The overall goal is to find the unitary $Q\in\U(n)$ that minimizes $\|\Log(Q^* \A)\|^2$ and $\norm{\sym\Log (Q^* \A)}^2$ over all possible logarithms.
Due to the non-uniqueness of the logarithm, we give the following as the formal statement of the minimization problem: 
\begin{align}
\min_{Q\in\U(n)}\norm{\Log (Q^* \A)}^2&:=\min_{Q\in\U(n)}\{  \norm{X}^2\in\R\,|\, \exp X=Q^* \A\}\, ,\label{min_defi0}\\
   \min_{Q\in\U(n)}\norm{\sym\Log (Q^* \A)}^2&:=\min_{Q\in\U(n)}\{  \norm{\sym X}^2\in\R\,|\, \exp X=Q^* \A\}\, \label{min_defi1}.
\end{align}
Our main result is the following. 
\begin{theorem}\label{thm:main}
 Let $Z\in \mathbb{C}^{n\times n}$ be a nonsingular matrix and let $Z=U_pH$ be its polar decomposition. Then 
\[
  \min_{Q\in U(n)} \norm{\Log Q^*Z}
=\min_{Q\in U(n)} \norm{\sym\Log Q^*Z}
=\norm{\log U_p^*Z}=\norm{\log H},
\]
for any $n$ when the norm is taken to be the spectral norm, and 
for $n\leq 3$ in the Frobenius norm. 
\end{theorem}

Our optimality result relies crucially on unitary invariance and a Bernstein-type trace inequality \cite{Bernstein88}
\begin{align}
  \tr{\exp X\, \exp{X^*}}\le \tr{\exp{(X+X^*)}}\, ,
\end{align}
for the matrix exponential. Together, these imply some algebraic conditions on the eigenvalues in case of the Frobenius matrix norm, which we exploit using 
a new sum of squared logarithms inequality \cite{Neff_log_inequality13}. For the spectral norm the analysis is considerably easier. 

This paper is organized as follows. In the remainder of this section we describe an application that motivated this work. In Section \ref{sec:complexscalar} we present two-dimensional analogues to our minimization problems in both, complex and real matrix representations, to illustrate the general approach and notation. In Section~\ref{sec:prep} we collect properties of the matrix logarithm and its Hermitian part. Section~\ref{sec:minimize} contains the main results where we discuss the unitary minimization \eqref{problem2}. From the complex case we then infer the real case in Section~\ref{sec:real} and finally discuss uniqueness in Section~\ref{sec:unique}. 

\emph{Notation}. 
$\sigma_i(X)=\sqrt{\lambda_i(X^*X)}$ denotes the $i$-th largest singular value of $X$. 
$\|X\|_2=\sigma_1(X)$ is the spectral matrix norm, $\norm{X}_F=\sqrt{\sum_{i,j=1}^n \abs{X_{ij}}^2}$ is the Frobenius matrix norm  with associated inner product $\Mprod{X}{Y}=\tr{X^* Y}$. 
The symbol $\id$ denotes the identity matrix.
An identity involving $\|\cdot\|$ without subscripts holds for any unitarily invariant norm. 
To avoid confusion between the unitary polar factor
and the singular value decomposition (SVD) of $\A=U\Sigma V^\ast$, 
$U_p$ with the subscript $p$ always denotes the unitary polar factor, while $U$
denotes the matrix of left singular vectors. Hence for example
$\A=U_{p}H=U\Sigma V^\ast $. $\U(n)$, $\Oo(n)$, $\GL(n,\C)$, $\GL^+(n,\R)$, $\SL(n)$ and $\SO(n)$ denote the group of complex unitary matrices, real orthogonal matrices, invertible complex matrices, invertible real matrices with positive determinant, the special linear group and the special orthogonal group, respectively. The set $\so(n)$ is the Lie-algebra of all $n\times n$ skew-symmetric matrices and $\sL(n)$ denotes the Lie-algebra of all $n\times n$ traceless matrices. The set of all $n\times n$ Hermitian matrices is $\mathbb{H}(n)$ and positive definite Hermitian matrices are denoted by $\mathbb{P}(n)$. We let  $\sym X=\frac{1}{2}(X^*+X)$ denote the Hermitian part of $X$ and $\skew X=\frac{1}{2}(X-X^*)$ the skew-Hermitian part of $X$ such that $X=\sym X+\skew X$. In general, $\Log Z$ with capital letter denotes any solution to $\exp X=Z$, while $\log Z$ denotes the principal logarithm.

\subsection{Application and practical motivation for the matrix logarithm}
 In this subsection we describe how our minimization problem concerning the matrix logarithm arises from new concepts in nonlinear elasticity theory and may find applications in generalized Procrustes problems. Readers interested only in the optimality result may continue reading Section \ref{sec:complexscalar}.

\subsubsection{Strain measures in linear and nonlinear elasticity}
Define the Euclidean distance ${\rm{dist}}_{\rm{euclid}}^2(X,Y):=\norm{X-Y}_F^2$, which is the length of the line segment joining $X$ and $Y$ in $\R^{n^2}$. We consider an elastic body which in a reference configuration occupies the bounded domain $\Omega\subset\R^3$. Deformations of the body are prescribed by mappings
\begin{align}
    \varphi:\Omega\mapsto \R^3\, ,
\end{align}
where $\varphi(x)$ denotes the deformed position of the material point $x\in\Omega$. Central to elasticity theory is the notion of strain. Strain is a measure of deformation such that no strain means that the body $\Omega$ has been moved rigidly in space. In linearized elasticity, one considers $\varphi(x)=x+u(x)$, where $u:\Omega\subset\R^3\mapsto \R^3$ is the displacement. The classical linearized strain measure is $\varepsilon:=\rsym\nabla u$. It appears through a matrix nearness problem
\begin{align}
    \dist_{\rm euclid}^2(\nabla u,\so(3)):=\min_{W\in\so(3)} \norm{\nabla u-W}_F^2=\norm{\rsym \nabla u}_F^2\, .
\end{align}
Indeed, $\rsym\nabla u$ qualifies as a linearized strain measure: if $\dist_{\rm euclid}^2(\nabla u,\so(3))=0$ then $u$ is a linearized rigid displacement of the form $u(x) = \widehat{W}x + \widehat{b}$ with fixed $\widehat{W}\in\so(3)$ and $\widehat{b}\in\R^3$.
This is the case since 
\begin{align}
      \dist_{\rm euclid}^2(\nabla u(x),\so(3))=0 \, \Rightarrow \, \nabla u(x)=W(x)\in\so(3)
      \end{align}
and $0=\Curl \nabla u(x)=\Curl W(x)$ implies that $W(x)$ is constant, see \cite{Neff_curl06}. 

In nonlinear elasticity theory one assumes that $\nabla \varphi\in\GL^+(3,\R)$ (no self-interpenetration of matter) and considers the matrix nearness problem
\begin{align}
    {\rm{dist}}_{\rm{euclid}}^2(\nabla \varphi,\SO(3)):=\min_{Q\in\SO(3)} \norm{\nabla \varphi-Q}_F^2=\min_{Q\in\SO(3)}\norm{Q^T\nabla\varphi-\id}_F^2\, .
\end{align}
From \eqref{unitary_polar_optimal} it  immediately follows that 
\begin{align}
\label{dist_euclid}
{\rm{dist}}_{\rm{euclid}}^2(\nabla\varphi,\SO(3) )=\norm{\sqrt{\nabla \varphi^T \nabla\varphi}-\id}_F^2\, . 
\end{align}
The term $\sqrt{\nabla\varphi^T \nabla\varphi}$ is called the right stretch tensor and $\sqrt{\nabla\varphi^T \nabla\varphi}-\id$ is called the Biot strain tensor. 
Indeed, the quantity $\sqrt{\nabla \varphi^T \nabla\varphi}-\id$ qualifies as a nonlinear strain measure: if $ {\rm{dist}}_{\rm{euclid}}^2(\nabla\varphi,\SO(3) )          =0$ then $\varphi$ is a rigid movement of the form $\varphi(x)=\widehat{Q}\,x+\widehat{b}$ with fixed $\widehat{Q}\in\SO(3)$ and $\widehat{b}\in\R^3$.
This is the case since 
\begin{align}
     {\rm{dist}}_{\rm{euclid}}^2(\nabla\varphi,\SO(3) )=0 \, \Rightarrow \, \nabla\varphi(x)=Q(x)\in\SO(3)
      \end{align}
and $0=\Curl \nabla \varphi(x)=\Curl Q(x)$ implies that $Q(x)$ is constant, see \cite{Neff_curl06}. 
Many other expressions can serve as strain measures. One classical example is the Hill-family \cite{Hill68,Hill70,Ogden70} of strain measures
\begin{align}
\label{strain_measure_family}
   a_m(\nabla\varphi):=  \begin{cases}
         \frac{1}{m}\left( \sqrt{\nabla \varphi^T \nabla\varphi}^{\,m}-\id\right)\, ,  & m\neq 0\\
         \log \sqrt{\nabla \varphi^T \nabla\varphi}\, , & m=0\, .
     \end{cases}    
\end{align}
The case $m=0$ is known as Hencky's strain measure \cite{Hencky28}. Note that the Taylor expansion $a_m(\id+\nabla u)=\sym\nabla u+O(u^2)$ coincide in the first-order approximation for all $m\in \R$.

In case of isotropic elasticity the formulation of a so-called boundary value problem of place may be based on postulating an elastic energy by integrating an $\SO(3)$-bi-invariant (isotropic and frame-indifferent) function $W:\R^{3\times3}\mapsto \R$ of the strain measure $a_m$ over $\Omega$
\begin{align}
    \mathcal{E}(\varphi):= \int_\Omega W(a_m(\nabla\varphi)) \,{\rm dx}\,,  \quad \varphi(x)_{\Gamma_D}=\varphi_0(x)\, 
\end{align}
and prescribing the boundary deformation $\varphi_0$ on the Dirichlet part $\Gamma_D\subset\partial\Omega$. The goal is to minimize $\mathcal{E}(\varphi)$ in a class of admissible functions. For example, choosing $m=1$ 
and $W(a_m)=\mu\, \norm{a_m}_F^2+\frac{\lambda}{2}\,  (\tr{a_m})^2$ leads to the isotropic Biot strain energy \cite{Neff_Biot07}
\begin{align}
     \int_\Omega \mu\, \norm{\sqrt{\nabla\varphi^T\nabla\varphi}-\id}_F^2+\frac{\lambda}{2}\left( \tr{\sqrt{\nabla\varphi^T\nabla\varphi}-\id}\right)^2\,{\rm dx}\
\end{align}
with Lam\'e constants $\mu,\lambda$. The corresponding Euler-Lagrange equations constitute a nonlinear, second order system of partial 
differential equations. For reasonable physical response of an elastic material Hill \cite{Hill68,Hill70,Ogden70} has argued that $W$ should be a convex function of the logarithmic strain measure $a_0(\nabla\varphi) =\log \sqrt{\nabla \varphi^T \nabla\varphi}$. This is the content of \emph{Hill's inequality}. 
Direct calculation shows that $a_0$ is the only strain measure among the family \eqref{strain_measure_family} that has the \emph{tension-compression symmetry}, i.e., for all unitarily invariant norms 
\begin{align}
    \norm{a_0(\nabla\varphi(x)^{-1})}=\norm{a_0(\nabla\varphi(x))}\, .
\end{align}

In his Ph.D. thesis \cite{Neff_Diss00} the first author was the first to observe that energies convex in the logarithmic strain measure $a_0(\nabla\varphi)$ are, in general, not rank-one convex. However, rank-one convexity is true in a large neighborhood of the identity \cite{Bruhns01}.

Assume for simplicity that we deal with an elastic material that can only sustain volume preserving deformations. Locally, we must have $\Det{\nabla\varphi(x)}=1$. Thus, for the deformation gradient $\nabla \varphi(x)\in\SL(3)$. On $\SL(3)$ the straight line $X+t(Y-X)$ joining $X,Y\in\SL(3)$ leaves the group. Thus, the Euclidean distance ${\rm{dist}}_{\rm{euclid} }^2(\nabla\varphi,\SO(3) )$ does not respect the group structure of $\SL(3)$.

Since the Euclidean distance \eqref{dist_euclid} is an arbitrary choice, novel approaches in nonlinear elasticity theory aim at putting more geometry (i.e. respecting the group structure of the deformation mappings) into the description of the strain a material endures.  In this context, it is natural to consider the strain measures induced by the geodesic distances stemming from choices for the Riemannian structure respecting also the algebraic group structure, which we introduce next.

\subsubsection{Geodesic distances}
In a connected Riemannian manifold $\mathcal{M}$ with Riemannian metric $g$, the length of a continuously differentiable curve $\gamma : [a,b] \mapsto \mathcal{M}$ is defined by
\begin{align}
   L(\gamma):&=\int\nolimits_a^b\sqrt{g_{\gamma(t)}(\dot{\gamma}(s),\dot{\gamma}(s))}\ds\, .
\end{align}
At every $X\in  \mathcal{M}$ the metric $g_X: T_X \mathcal{M}\times T_X\mathcal{M} \mapsto\R$ is a positive definite, symmetric bilinear form on the tangent space $T_X \mathcal{M}$. The distance $ \dist_{{\rm geod},\mathcal{M}}(X,Y)$ between two points $X$ and $Y$ of $\mathcal{M}$ is defined as the infimum of the length taken over all continuous, piecewise continuously differentiable curves $\gamma : [a,b] \mapsto \mathcal{M}$  such that $\gamma(a) = X$ and $\gamma(b) = Y$. 
See~\cite{andruchow2011left} for more discussion on the geodesics distance.
With this definition of distance, geodesics in a Riemannian manifold are the locally distance-minimizing paths, in the above sense. Regarding $\mathcal{M}=\SL(3)$ as a Riemannian manifold equipped with the metric associated to one of the positive definite quadratic forms of the family
\begin{align}
\label{Riemannian_metric}
  g_X(\xi,\xi):&=\mu\, \norm{{\rm sym} (X^{-1}\xi)}_F^2+\mu_c\, \norm{{\rm skew}(X^{-1} \xi)}_F^2\, ,\; \xi\in T_X\SL(3) \, 
\end{align}
for  all $\mu,\mu_c>0$,
where we drop the subscript $_\ast$ in $\rm sym_\ast$ when the matrix is real,
we have 
$\gamma^{-1}(t)\dot \gamma(t)\in T_\id\SL(3)=\sL(3)$
 (by direct calculation, $\sL(3)$ denotes the trace free $\R^{3\times 3}$-matrices) and 
\begin{align}  
  g_{\gamma(t)}(\dot{\gamma}(t),\dot{\gamma}(t))&=\mu\, \norm{{\rm sym} (\gamma^{-1}(t)\dot{\gamma}(t))}_F^2+\mu_c\, \norm{{\rm skew}(\gamma^{-1}(t) \dot{\gamma}(t))}_F^2\,  .    
\end{align}
It is clear that
\begin{align}
   \forall\; \mu,\mu_c>0:\quad    \mu\, \norm{{\rm sym} \,Y}_F^2+\mu_c\, \norm{{\rm skew} \,Y}_F^2
     \end{align}
is a norm on $\sL(3)$. For such a choice of metric we then obtain an associated Riemannian distance metric
\begin{align}
\label{definition_geodesic}
    \dist_{{\rm geod},\SL(3)}(X,Y)=\inf\{ L(\gamma),\; \gamma(a)=X\, ,\; \gamma(b)=Y\, \}\, .
\end{align}
This construction ensures the validity of the triangle inequality \cite[p.14]{Jost02b}. The geodesics on $\SL(3)$ for the family of metrics \eqref{Riemannian_metric} have been computed in \cite{Mielke02b} in the context of dissipation distances in elasto-plasticity.

With this preparation, it is now natural to consider the strain measure induced by the geodesic distance. For a given deformation gradient $\nabla\varphi\in\SL(3)$ we thus compute the distance to the nearest orthogonal matrix in the geodesic distance \eqref{definition_geodesic} on the Riemanian manifold and matrix Lie-group $\SL(3)$, i.e.,
\begin{align}
\label{dist_geod_general}
     \dist_{{\rm geod},\SL(3)}^2(\nabla\varphi,\SO(3)):=\min_{Q\in\SO(3)}    \dist_{{\rm geod},\SL(3)}^2(\nabla\varphi,Q)    \, . 
 \end{align}
 It is clear that this defines a strain measure, since $\dist_{{\rm geod},\SL(3)}^2(\nabla\varphi(x),\SO(3))=0$ implies $\nabla\varphi(x)\in\SO(3)$, whence $\varphi(x)=\widehat{Q}x+\widehat{b}$. Fortunately, the minimization on the right hand side in \eqref{dist_geod_general} can be carried out although the explicit distances  $\dist_{{\rm geod},\SL(3)}^2(\nabla\varphi,Q)$ for given $Q\in\SO(3)$ remain unknown to us. In \cite{Neff_Osterbrink_hencky13} it is shown that 
 \begin{align}
\label{dist_geod}
     \min_{Q\in\SO(3)}    \dist_{{\rm geod},\SL(3)}^2(\nabla\varphi,Q) =\min_{Q\in\SO(3)} \| \Log(Q^T \nabla \varphi)\|_F^2\, .
 \end{align}
Recall that $\Log Z$ denotes any matrix logarithm, one of the many solutions $X$ to $\exp X = Z$. By contrast, $\log Z$ denotes the principal logarithm, see Section~\ref{sec:princ}.
The last equality constitutes the basic motivation for this work, where we solve the minimization problem on the right hand side of \eqref{dist_geod} and determine thus the precise form of the geodesic strain measure. As a result of this paper it turns out that 
\begin{align}
         \dist_{{\rm geod},\SL(3)}^2(\nabla\varphi,\SO(3))=\norm{\log\sqrt{\nabla\varphi^T\nabla\varphi} }_F^2\, ,
\end{align}
which is nothing else but a quadratic expression in Hencky's strain measure \eqref{strain_measure_family} and therefore satisfying Hill's inequality.\\

Geodesic distance measures have appeared recently in many other applications: for example, one considers a geodesic distance on the Riemannian manifold of the cone of positive definite matrices $\mathbb{P}(n)$ (which is a Lie-group but not w.r.t. the usual matrix multiplication) \cite{Bhatia06,Moakher05} given by%
\begin{align}
\label{dist_geod_pos}
 \dist_{{\rm geod},\mathbb{P}(n)}^2(P_1,P_2):=\norm{\log(P_1^{-1/2} P_2 P_1^{-1/2})}_F^2          \, .
 \end{align}
 Another distance, the so-called log-Euclidean metric on $\mathbb{P}(n)$ 
 \begin{align}
 \label{log_Euclidean_metric}
 \dist_{\log,\rm{euclid},\mathbb{P}(n)}^2&(P_1,P_2):=\norm{\log P_2-\log P_1}_F^2\notag\\
      &   \left( \text{in general}\neq  \norm{\log (P_1^{-1} P_2)}_F^2=\dist_{\log,\rm{euclid},\mathbb{P}(n)}^2(P_1^{-1}P_2,\id)  \right)
  \end{align}
 is proposed in \cite{Arsigny06}. Both formulas find application in diffusion tensor imaging or in fitting of positive definite elasticity tensors. The geodesic distance on the compact matrix Lie-group $\SO(n)$ is also well known, and it has important applications in the interpolation and filtering of experimental data given on $\SO(3)$, see e.g. \cite{Moakher02}
\begin{align}
  \label{dist_geod_SO}
      \dist_{{\rm geod},\SO(n)}^2(Q_1,Q_2):=\norm{\log (Q_1^{-1} Q_2)}_F^2\, , \quad -1\not\in {\rm spec}(Q_1^{-1} Q_2)\, .
 \end{align}     
Here ${\rm spec}(X)$ denotes the set of eigenvalues of the matrix $X$. 
In cases \eqref{dist_geod_pos}, \eqref{log_Euclidean_metric}, \eqref{dist_geod_SO} it is, contrary to \eqref{dist_geod}, the principal matrix logarithm that appears naturally. A common and desirable feature of all distance measures involving the logarithm presented above, setting them apart from the Euclidean distance, is invariance under inversion: ${\rm d}(X,\id)={\rm d}(X^{-1},\id)$ and ${\rm d}(X,0)=+\infty$. We note in passing that 
\begin{align}
    \label{logarithmic_one_parameter_distance}
      {\rm d}_{\log,\GL^+(n,\R)}^2(X,Y):=\norm{\Log (X^{-1} Y)}_F^2
      \end{align}
does not satisfy the triangle inequality and thus it cannot be a Riemannian distance metric on $\GL^+(n,\R)$. Further, $X^{-1}Y$ is in general not in the domain of definition of the principal matrix logarithm. If applicable, the expression \eqref{logarithmic_one_parameter_distance} measures in fact the length of curves $\gamma:[0,1]\mapsto \mathcal{M},\gamma(0)=X, \gamma(1)=Y$ defining one-parameter groups $\gamma(s)=X\, \exp(s \, \Log(X^{-1}Y) )$ on the matrix Lie-group $\mathcal{M}$. Note that it is only if the manifold $\mathcal{M}$ is a compact matrix Lie-group (like e.g. $\SO(n)$) equipped with a bi-invariant Riemannian metric that the geodesics are precisely one-parameter subgroups \cite[Prop.9]{Neill83}. This point is sometimes overlooked in the literature.

\subsubsection{A geodesic orthogonal Procrustes problem on $\SL(3)$}
The Euclidean orthogonal Procrustes problem for $\A,B\in\SL(3)$
\begin{align}
   \min_{Q\in O(3)} \dist_{\rm euclid}^2(\A,B Q)=\min_{Q\in O(3)} \norm{\A-B Q}_F^2
\end{align}
has as solution the unitary polar factor of $B^*\A$ \cite[Ch.~12]{Golubbookori}. However, any linear transformation of $\A$ and $B$ will yield another optimal unitary matrix. This deficiency can be circumvented by considering the straightforward extension to the geodesic case
\begin{align}
   \min_{Q\in O(3)} \dist_{{\rm geod},\SL(3)}^2(\A,B Q)\, .
\end{align}
In contrast to the Euclidean distance, the geodesic distance is by construction $\SL(3)$-left-invariant:
\begin{align}
   \dist_{{\rm geod},\SL(3)}^2(X,Y)=\dist_{{\rm geod},\SL(3)}^2(B X,B Y)\quad   \mbox{for\ all\ }
\, B\in \SL(3),
  \end{align}   
and therefore we have 
\begin{align}
   \min_{Q\in O(3)} \dist_{{\rm geod},\SL(3)}^2(\A,B Q)= \min_{Q\in O(3)} \dist_{{\rm geod},\SL(3)}^2(B^{-1}\A,Q) 
\end{align}
with ``another" geodesic optimal solution: the unitary polar factor of $B^{-1}\A$, according to the results of this paper. A more detailed description of this additional optimality result as well as its application towards elasticity theory can be found in \cite{Neff_Eidel_Osterbrink_2013}.

\section{Prelude on optimal rotations in the complex plane}
\label{sec:complexscalar}
Let us turn to the optimal rotation problem, the first term of \eqref{problem2}:
\begin{align}
 \min_{Q\in \U(n)} \| \Log(Q^* \A)\|^2\, .
 \end{align}
 In order to get hands on this problem we first consider the scalar case. 
It serves as a useful preparation for the matrix case, as we follow the same logical sequence in the next section.
 We may always identify the punctured complex plane $\C\setminus\{0\}=:\C^\times=\GL(1,\C)$ 
 with the two-dimensional conformal special orthogonal group $\CO(2)\subset \GL^+(2,\R)$ through the mapping
\begin{align}
\label{matrix_complex}
 z=a+i\, b\quad \mapsto \quad 
 \A\in \CO(2):=\{ \begin{bmatrix}
     a & b\\
     -b & a
 \end{bmatrix}, \quad a^2+b^2\neq 0 \} \, .
 \end{align}
Let us define a norm $\norm{\cdot}_{\CO}$ on $\CO(2)$. We set $\norm{X}_{\CO}^2:=\frac{1}{2}\norm{X}_{F}^2=\frac{1}{2}\tr{X^TX}$. 

Next we introduce the logarithm. For every invertible $z\in\C\setminus\{0\}=:\C^\times$ there always exists a solution to $e^\eta=z$ and we call $\eta\in\C$ the natural complex logarithm $\Log_\C(z)$ of $z$. However, this logarithm may not be unique, depending on the unwinding number \cite[p.269]{Higham2008}. The definition of the natural logarithm has some well known deficiencies: the formula $\Log_\C(w^z)=z\Log_\C(w)$ does not hold, since, e.g. $i\, \pi=\Log_\C(-1)=\Log_\C((-i)^2)\neq 2 \Log_\C(-i)=2(\frac{-i\, \pi}{2})=-i\,\pi$. Therefore the principal complex logarithm \cite[p.79]{Bernstein2009} 
\begin{align}
  \logp: \C^\times \to \{\, z\in\C\, |\, -\pi<\Im z\le  \pi\, \}
\end{align}
is defined as the unique solution $\eta\in \C$ of 
\begin{align}
  e^\eta=z\quad \Leftrightarrow\quad  \eta=\logp(z):=\log\abs{z}+i\, \arg(z)\, ,
\end{align}
such that the argument $\arg(z)\in (-\pi,\pi]$.%
\footnote{For example $\logp(-1)=i\, \pi$ since $e^{i\pi}=-1$. Otherwise, the complex logarithm always exists but may not be unique, e.g. $e^{-i\pi}=\frac{1}{e^{i\, \pi}}=-1$. Hence $\Log_\C(-1)=\{i\, \pi, -i\, \pi, \ldots\}$. 
For scalars, our definition of the principal complex logarithm can be applied to negative real arguments. However, in the matrix setting the principal matrix logarithm is defined only for invertible matrices which do not have negative real eigenvalues.}
The principal complex logarithm is continuous (indeed holomorphic) only on the smaller set $\C\setminus (-\infty,0] $. Let us define the set $\mathcal{D}:=\{z\in\C\, |\, \abs{z-1}<1\,\}$. In order to avoid unnecessary complications at this point, we introduce a further open set, the ``near identity subset''
\[
	\mathcal{D}^\sharp := \{z\in\C\, |\, \abs{z-1}<\sqrt{2}-1\,\} \subset \mathcal{D}\,,
\]
which is defined such that $1\in\mathcal{D}^\sharp$ and $z_1, z_2\in\mathcal{D}^\sharp$ implies $z_1 z_2\in\mathcal{D}$ and $z_1^{-1}\in\mathcal{D}$.
On $\mathcal{D}^\sharp\subset\C^\times$ all the usual rules for the logarithm apply.
Simplifying further, on $\R^+\setminus\{0\}$ \emph{all} the logarithmic distance measures encountered in the introduction coincide with the logarithmic metric \cite[p.109]{Tarantola06} (the "hyperbolic distance" \cite[p.735]{Moakher05})
\begin{align}\label{real_logarithmic_metric}
 \dist_{\log,\R^+}^2(x,y):&=\abs{\log (x^{-1} y)}^2=\abs{\log y-\log x}^2, \ \,\notag\\
 \dist_{\log,\R^+}^2(x,1)&=\abs{\log \abs{x}}^2.
 \end{align}
This metric can still be extended to a metric on $\mathcal{D}^\sharp$ through
\begin{align}
\label{d_sharp_logarithmic_metric_definition}
 \dist_{\log,\mathcal{D}^\sharp}^2(z_1,z_2):&=\abs{\logp (z_1^{-1} z_2) }^2\, , \;z_1, z_2\in\mathcal{D}^\sharp\, ,\\
  \dist_{\log,\mathcal{D}^\sharp}^2(r_1 e^{i\vartheta_1},r_2 e^{i\vartheta_2})&= \abs{ \log (r_1^{-1} r_2)}^2+\abs{\vartheta_1-\vartheta_2}^2\, ,
  \quad r_1 e^{i\vartheta_1}, r_2 e^{i\vartheta_2}\in\mathcal{D}^\sharp\, .\notag
    \end{align}
Further, for $z\in\mathcal{D}^\sharp$ we find a formula similar to \eqref{real_logarithmic_metric} by taking the distance of $z$ to the set of all $y\in\mathcal{D}^\sharp$ with $\abs{y}=1$ instead of the distance to $1$:
\[
     \min_{e^{i\vartheta}\in\mathcal{D}^\sharp}\dist_{\log,\mathcal{D}^\sharp}^2(z,e^{i\vartheta})=\abs{ \log \abs{z}}^2\, .
\]
  We remark, however, that we cannot simply extend \eqref{d_sharp_logarithmic_metric_definition} to a metric on $\C^\times$ due to the periodicity of the complex exponential.
Let us also define a log-Euclidean distance metric on $\C^\times$, continuous only on $\C\setminus (-\infty,0]$, in analogy with \eqref{log_Euclidean_metric}
   \begin{align}
       \dist_{\log,\rm{euclid},\C^\times}^2(z_1,z_2):&=\abs{\logp z_2-\logp z_1}^2
        =\left|\log \frac{\abs{z_2}}{\abs{z_1}}\right|^2+\abs{\arg(z_2)-\arg(z_1)}^2.     
       \end{align}   
    The identity 
  \begin{align}  
    \label{complex_commutator_identity}
    \dist_{\log,\mathcal{D}^\sharp}^2(z_1,z_2)=  \dist_{\log,\rm{euclid},\C^\times}^2(z_1,z_2)
  \end{align}  
on $\mathcal{D}^\sharp$ is obvious, although it is not well-posed on $\C^\times$. With this preparation, we now approach our minimization problem in terms of $\CO(2)$ versus $\C^\times$. For given $\A\in\CO(2)$ we find that the following minimization problems are equivalent, meaning that if we identify $Q\in\SO(2)$ with the corresponding complex number $e^{i\vartheta}\in\C^\times$ using \eqref{matrix_complex}, the minimizing arguments are equal:
\begin{align}
\label{equivalence_matrix_complex}
    \min_{Q\in\SO(2)}\norm{\Log(Q^T \A)}_{\CO}^2 \;\sim\;
\min_{e^{i\vartheta}\in\C^\times}\abs{\Log_{\C}(e^{-i\vartheta} z)}^2\,.
\end{align}
Here $\Log_{\C}$ is as defined below in \eqref{complex_multivalued_definition}.
It is important to avoid the additive representation inherent in $ \dist_{\log,\rm{euclid},\C^\times}^2(z_1,z_2)$, 
because in the general matrix setting $Q$ and $\A$ will in general not commute and the equivalence of the problems in \eqref{equivalence_matrix_complex} is then lost.
%

In order to give the minimization problem \eqref{equivalence_matrix_complex} a precise sense, we define
\begin{equation}
\label{complex_multivalued_definition}
 \min_{\vartheta\in(-\pi,\pi]}\abs{\Log_{\C}(e^{-i\vartheta} z)}^2 := \min_{\vartheta\in(-\pi,\pi]} \min\,\{ \abs{w}^2\, | \, e^w=e^{-i\vartheta} z\, \}
\end{equation}
as the minimum over \emph{all} logarithms of $e^{-i\vartheta} z$.
Dropping the second ``$\min$'' for better readability we find
\begin{align}
\label{complex_multivalued_alternative_definition}
 &\min_{\vartheta\in(-\pi,\pi]}\abs{\Log_{\C}(e^{-i\vartheta} z)}^2\\
 &\qquad\quad= \min_{\vartheta\in(-\pi,\pi]}\{ \abs{w}^2\, | \, e^w=e^{-i\vartheta} z\, \}=
  \min_{\vartheta\in(-\pi,\pi]}\{ \abs{w}^2\, | \, e^w=e^{-i\vartheta} e^{i\arg(z)}\abs{z}\, \} \notag\\
  &\qquad\quad=\min_{\tilde{\vartheta}\in(-\pi,\pi]}\{ \abs{w}^2\, | \, e^w=e^{-i\tilde{\vartheta}} \abs{z}\, \} 
  = \min_{\vartheta\in(-\pi,\pi]}\abs{\Log_{\C}(e^{-i\vartheta} \abs{z})}^2  \, . 
    \end{align}
The solution of this minimization problem is again
  $\abs{\log\abs{z}}^2$, since
 $\min_{\vartheta\in(-\pi,\pi]}\abs{\Log_{\C}(e^{-i\vartheta} \abs{z})}^2=\min_{\vartheta\in(-\pi,\pi]}\abs{\log\abs{z}+i(-\vartheta)}^2=  \min_{\vartheta\in(-\pi,\pi]}\abs{\log\abs{z}}^2+\abs{\vartheta}^2  
  =\abs{\log\abs{z}}^2$. 
 However, 
our goal is to introduce an argument that can be generalized to the non-commutative matrix setting. From $\abs{z}\ge \abs{\Re(z)}$ it follows that
\begin{align}
\label{minimum_complex}
   \min_{\vartheta\in(-\pi,\pi]}\abs{\Log_{\C}(e^{-i\vartheta} z)}^2&\ge  \min_{\vartheta\in(-\pi,\pi]}\abs{\Re(\Log_{\C}(e^{-i\vartheta} z))}^2 =\abs{\log\abs{z}} ^2 \, ,
   \end{align}
where we used the result \eqref{minimum_complex_symmetric} below for the last equality.
  The minimum for $\vartheta\in(-\pi,\pi]$ is achieved if and only if $\vartheta=\arg(z)$ since $\arg(z)\in (-\pi,\pi]$ and we are looking only for $\vartheta\in(-\pi,\pi]$. Thus 
     \begin{align}
    \min_{\vartheta\in(-\pi,\pi]}\abs{\Log_{\C}(e^{-i\vartheta} z)}^2= \abs{\log \abs{z}} ^2 \, .
 \end{align}   
 The unique optimal rotation $Q(\vartheta)\in\SO(2)$ is given by the polar factor $U_p$ through $\vartheta=\arg(z)$ and the minimum is $\abs{\log\abs{z}} ^2$, which corresponds to $\min_{Q\in\SO(2)}\norm{\Log(Q^T \A)}_{\CO}^2=\norm{\log \sqrt{\A^T \A}}_{\CO}^2$. 

Next, consider the symmetric minimization problem in \eqref{problem2} for given $\A\in\CO(2)$ and its equivalent representation in $\C^\times$:
\begin{align}
\label{complex_real_matrix_sym_equaivanlence}
    \min_{Q\in\SO(2)}\norm{\sym\Log(Q^T \A)}_{\CO}^2 \quad \sim\quad  
   \min_{\vartheta\in(-\pi,\pi]}\abs{\Re(\Log_{\C}(e^{-i\vartheta} z))}^2   \, .
\end{align}
Note that the expression $ \dist_{\log,\Re,\C^\times}(z_1,z_2):=\abs{\Re \Log_{\C} (z_1^{-1} z_2) }$ does not define a metric, even when restricted to $\mathcal{D}^\sharp$. As before, we define
\begin{align}
  \label{complex_multivalued_definition_real}
   \min_{\vartheta\in(-\pi,\pi]}\abs{\Re(\Log_{\C}(e^{-i\vartheta} z))}^2:=\min_{\vartheta\in(-\pi,\pi]}\{ \abs{\Re w}^2\, | \, e^w=e^{-i\vartheta} z\, \} \,  \end{align}
and obtain
\begin{align}
\label{minimum_complex_symmetric}
   \min_{\vartheta\in(-\pi,\pi]}\abs{\Re(\Log_{\C}(e^{-i\vartheta} z))}^2
   &= \min_{\vartheta\in(-\pi,\pi]}\abs{\Re(\Log_{\C}(e^{-i\vartheta} \abs{z} e^{i\, \arg(z)}))}^2 \notag\\
  &=\min_{\vartheta\in(-\pi,\pi]}\abs{\Re(\Log_{\C}(\abs{z} e^{i\, (\arg(z)-\vartheta)}))}^2 \\
  &=\min_{\vartheta\in(-\pi,\pi]} \{ \abs{\Re(\log\abs{z}+i ( \arg(z)-\vartheta+2\pi\, k))}^2 \, , k\in\N\} \notag\\
  &=\abs{\log\abs{z}} ^2 \, .\notag
  \end{align}
Thus the minimum is again realized by the polar factor $U_p$, but note that the optimal rotation is completely undetermined, since $\vartheta$ is not constrained in the problem. Despite the logarithm $\Log_\C$ being multivalued, this formulation of the minimization problem circumvents the problem of the branch points of the natural complex logarithm. This observation 
 suggests that considering the generalization of \eqref{minimum_complex_symmetric}, i.e. $\min_{Q\in\U(n)}\norm{\sym\Log Q^* \A}^2$  in the first place is helpful also for the general matrix problem. This is indeed the case.

With this preparation we now turn to the general, non-commutative matrix setting.

\section{Preparation for the general complex matrix setting}\label{sec:prep}

\subsection{Multivalued formulation}
For every nonsingular $Z\in\GL(n,\C)$ there exists a solution $X\in\C^{n\times n}$ to $\exp X=Z$ which we call a logarithm $X=\Log(\A)$ of $Z$. 
As for scalars, the matrix logarithm is multivalued depending on the unwinding number~\cite[p.~270]{Higham2008} since in general, a nonsingular real or complex matrix may have an infinite number of real or complex logarithms. The goal, nevertheless, 
is to find the unitary $Q\in\U(n)$ that minimizes $\|\Log(Q^* \A)\|^2$ and $\norm{\sym\Log (Q^* \A)}^2$ over all possible logarithms.

Since $\|\Log(Q^* \A)\|,\, \norm{\sym\Log (Q^* \A)}^2\ge 0$, it is clear that both infima exist.  Moreover, $\U(n)$ is compact and connected. One problematic aspect is that $\U(n)$ is a non-convex set and the function $X\mapsto\norm{\Log X}^2$ is non-convex. Since, in addition, the multivalued matrix logarithm may fail to be continuous, at this point we cannot even claim the existence of minimizers.


We first observe that without loss of generality we may assume that $\A\in \GL(n,\C)$ is real, diagonal and positive definite. 
To see this, consider the unique polar decomposition $\A=U_p \, H$ 
and the eigenvalue decomposition 
$H=V D V^*$ 
for real diagonal positive $D=\diag(d_1,\ldots, d_n)$. Then, in complete analogy to
\eqref{complex_multivalued_alternative_definition},
\begin{align}
 \min_{Q\in\U(n)}\norm{\sym\Log (Q^* \A)}^2&=\min_{Q\in\U(n)}  \{  \norm{\sym X}^2 \,|\, \exp X=Q^* \A\}\notag\\
 &=\min_{Q\in\U(n)} \{  \norm{\sym X}^2 \,|\, \exp X=Q^* U_p  H\}\notag\\
&=\min_{Q\in\U(n)} \{  \norm{\sym X}^2 \,|\, \exp X=Q^*  U_pV D V^*\}\notag\\
&=\min_{Q\in\U(n)} \{  \norm{\sym X}^2 \,|\, V^* (\exp X) V=V^* Q^* U_pV D\}\notag\\
&=\min_{Q\in\U(n)} \{  \norm{\sym X}^2 \,|\, \exp (V^*XV)=V^* Q^* U_pV D\}\notag\\
&=\min_{\widetilde{Q}\in\U(n)} \{  \norm{\sym X}^2 \,|\,  \exp (V^*XV)=\widetilde{Q}^* D \}\notag\\
&=\min_{\widetilde{Q}\in\U(n)} \{  \norm{V^*(\sym X)V}^2 \,|\,  \exp (V^*XV)=\widetilde{Q}^* D \}\notag\\
&=\min_{\widetilde{Q}\in\U(n)} \{  \norm{\sym (V^*X V)}^2 \,|\,  \exp (V^*XV)=\widetilde{Q}^* D \}\notag\\
&=\min_{\widetilde{Q}\in\U(n)}\{  \norm{\sym (\widetilde{X}) }^2 \,|\,  \exp (\widetilde{X})=\widetilde{Q}^* D \}\notag\\
&=\min_{Q\in\U(n)} \{  \norm{\sym X}^2 \,|\, \exp X=Q^* D\}\notag\\
&= \min_{Q\in\U(n)}\norm{\sym\Log Q^*D}^2\, ,\notag
 \end{align}
where we used the unitary invariance for any unitarily invariant matrix norm and the fact that $X\mapsto \sym X$ and $X\mapsto \exp X$ are isotropic functions, i.e. invariant under congruence with orthogonal/unitary transformations $f(V^* X V)=V^* f(X) V$ for all unitary $V$. If the minimum is achieved for $Q=\id$ in $\min_{Q\in \U(n)}\norm{\sym \Log (Q^* D) }^2$ then this corresponds to $Q=U_p$ in $\min_{Q\in \U(n)}\norm{\sym \Log Q^* \A}^2$. Therefore, in the following we assume  that $D=\diag(d_1,\ldots d_n)$ with $d_1\ge d_2\ge \ldots \ge d_n>0$.

\subsection{Some properties of the matrix exponential $\exp$ and matrix logarithm $\Log$}

Let $Q\in \U(n)$. Then the following equalities hold for all $X\in\C^{n\times n}$.
\begin{align}
    \exp(Q^* X Q)&=Q^* \exp(X) \, Q\, ,\quad\text{definition of $\exp$, \cite[p.715]{Bernstein2009}} ,\\
     Q^* \Log(X) \, Q&\quad \text{is a logarithm of}\quad Q^* X Q\,,\\
     \det{Q^* X Q}&=\det X\, \, ,\\
\exp (-X)&=\exp(X)^{-1}\, , \quad\quad \text{series definition of $\exp$, \cite[p.713]{Bernstein2009}}\, ,\notag\\ 
     \exp\Log X&=X\,,\quad \quad \quad\text{for any matrix logarithm}\, , \\
 \det{ \exp X}&=e^\tr{X}\, ,\quad\quad  \text{\cite[p.712]{Bernstein2009}}\, ,\\
\forall\, Y\in \C^{n\times n},\det{Y}&\neq 0: \, \det{Y}=e^{\tr{\Log Y}} \quad\text{for any matrix logarithm \cite{Higham2008}}\, .\notag 
\end{align}
A major difficulty in the multivalued matrix logarithm case arises from
\begin{align}
\forall\,X\in \C^{n\times n}:\quad \Log\exp X&\neq X  \quad\text{in general, without further assumptions}\, .
\end{align}

\subsection{Properties of the principal matrix-logarithm $\logp$}\label{sec:princ}
Let $X\in \C^{n\times n}$, and assume that $X$ has no real eigenvalues in $(-\infty,0]$. The principal matrix logarithm of $X$ is the unique logarithm of $X$ (the unique solution $Y\in\C^{n\times n}$ of $\exp Y=X$) whose eigenvalues are elements of the strip $\{z\in \C:\; -\pi< \Im(z) <\pi\}$. If $X\in\R^{n\times n}$ and $X$ has no eigenvalues on the closed negative real axis $\R^-=(-\infty,0]$, then the principal matrix logarithm is real. 
Recall that $\logp X$ is the principal logarithm and $\Log X$ denotes one of the many solutions to $\exp Y=X$. 

The following statements apply strictly only to the principal matrix logarithm \cite[p.721]{Bernstein2009}:
\begin{align}
\label{rules_principal_logarithm}
\logp\exp X&=X\quad \text{if and only if $\abs{\Im \lambda}<\pi$ for all $\lambda\in \rm{spec}(X)$}\, , \notag\\
\logp(X^{\alpha})&=\alpha \log X\, ,\quad\quad \alpha\in[-1,1]\, ,\notag\\
 \logp(Q^* X Q)&=Q^* \logp(X) \, Q\, ,\quad \forall\, Q\in\U(n)\, .
\end{align}
Let us define the set of Hermitian matrices $\mathbb{H}(n):=\{X\in\C^{n\times n}\, | \, X^*=X\, \}$ and the set $\mathbb{P}(n)$ of positive definite Hermitian matrices consisting of all Hermitian matrices with only positive eigenvalues. 
The mapping
\begin{align}
\label{exp_positive}
   \exp: \mathbb{H}(n)\mapsto \mathbb{P}(n)
\end{align}
is bijective \cite[p.719]{Bernstein2009}. 
In particular, 
$\Log \exp\sym X$ is uniquely defined for any $X\in\C^{n\times n}$ 
up to additions by multiples of $2\pi i$ to each eigenvalue
and any matrix logarithm and therefore we have
\begin{align}
\label{log_on_psym_1}
    \forall\, H\in\mathbb{H}(n):\quad \sym \Log H&=\logp H\, ,\notag\\
   \forall\,X\in \C^{n\times n}:\quad \logp \exp\sym X&=\sym X\, ,  \\
    \forall\,X\in \C^{n\times n}:\quad \sym\Log \exp\sym X&=\sym X\, .\notag
\end{align}
 Since $\exp\sym X$ is positive definite, it follows from \eqref{rules_principal_logarithm} also that
  \begin{align}    
  \label{log_on_psym_2}
 \forall\,X\in \C^{n\times n}:\quad Q^*(\logp\exp\sym X) Q&=\logp(Q^*(\exp\sym X) Q) \, .
\end{align}

\section{Minimizing $\| \Log(Q^* \A))\|^2$}\label{sec:minimize}
Our starting point is, in analogy with the complex case \eqref{complex_real_matrix_sym_equaivanlence}, the problem of minimizing 
\begin{align*}
\min_{Q\in \U(n)}\|\sym(\Log(Q^* \A))\|^2\, ,
\end{align*}
where 
$\sym(X)=(X^*+X)/2$ is the Hermitian part of $X$. As we will see, a solution of this problem will already imply the full statement, similar to the complex case, see \eqref{minimum_complex}. For every complex number $z$, we have 
\begin{align}
\label{complex_number_estimate}
     \abs{e^z}=e^{\Re z}=\abs{e^{\Re z}   }\le \abs{e^{\Re z}   }\, .
\end{align}
While the last inequality in \eqref{complex_number_estimate} is superfluous it is in fact the ``inequality'' $\abs{e^z}\le \abs{e^{\Re z}}$ that can be generalized to the matrix case. The key result is an inequality of Bhatia \cite[Thm.~IX.3.1]{Bhatia97}, 
\begin{equation}  
\label{bhatia}
\forall\, X\in\C^{n\times n}:\quad \|\exp X\|^2\leq   \|\exp \sym X\|^2
\end{equation}
for any unitarily invariant norm, cf. \cite[Thm.~10.11]{Higham2008}. The result \eqref{bhatia} is a generalization of Bernstein's trace inequality 
for the matrix exponential: in terms of the Frobenius matrix norm it holds
\begin{align*}
   \norm{\exp{X}}_F^2=\tr{\exp X\, \exp{X^*}}\le \tr{\exp{(X+X^*)}}=\norm{\exp{\sym X}}_F^2\, ,
\end{align*}
with equality if and only if $X$ is normal \cite[p.756]{Bernstein2009}, \cite[p.515]{Horn91}. For the case of the spectral norm the inequality 
\eqref{bhatia} is already given by Dahlquist \cite[(1.3.8)]{Dahlquist59}. We note that the well-known Golden-Thompson inequalities \cite[p.761]{Bernstein2009},\cite[Cor.6.5.22(3)]{Horn91}: 
\begin{align*}
\forall\; X,Y\in \mathbb{H}(n): \quad \tr{\exp(X+Y)}\le \tr{\exp(X)\,\exp(Y)}  
\end{align*}
seem (misleadingly) to suggest the reverse inequality.\\

Consider for the moment any unitarily invariant norm, any $Q\in \U(n)$, the positive real diagonal matrix $D$ as before and any 
matrix logarithm $\Log$. Then it holds
\begin{align}
\label{exp_estimate_1}
\norm{\exp(\sym\Log Q^* D)}^2&\ge \norm{\exp(\Log Q^* D)}^2=\norm{Q^* D}^2=\norm{D}^2 \, ,  
\end{align}
due to inequality \eqref{bhatia} and
\begin{align}
\label{exp_estimate_2}
\norm{\exp(-\sym\Log Q^* D)}^2&=\norm{\exp(\sym(-\Log Q^* D))}^2\notag\\
&\ge \norm{\exp((-\Log Q^* D))}^2
=\norm{(\exp(\Log Q^* D))^{-1}}^2\notag\\
&=\norm{(Q^* D)^{-1}}^2=\norm{D^{-1} (Q^*)^{-1}}^2=\norm{D^{-1}}^2   \, ,
\end{align}
where we used \eqref{bhatia} again. Note that we did not use $-\Log X=\Log(X^{-1})$ (which may be wrong, depending on the unwinding number). 

Moreover, we note that for any $Q\in \U(n)$ we have
\begin{align}
\label{det_equality}
0<\det{\exp(\sym\Log Q^* D)}&=e^\tr{\sym\Log Q^* D}=e^{\Re\tr{\Log Q^* D}}\notag\\
                                                     &=\abs{e^{\Re\tr{\Log Q^* D}}}=\abs{e^\tr{\Log Q^* D}}    \\
    &=\abs{\det{Q^* D}}=\abs{\det{Q^*}\det{D} }\notag\\
    &=\abs{\det{Q^*}}\,\abs{\det{D}}=\abs{\det D}=\det D\, , \notag
\end{align}
where we used the fact that 
\begin{align}
\label{log_det_computation}
e^{\tr{X}}&=\det{\exp{X}}\,, \quad X=\Log Q^*D\; \Rightarrow \; \notag\\
&e^{\tr{\Log Q^*D}}=\det{\exp{\Log Q^* D}}=\det{Q^*D}\, ,
\end{align}
is valid for any solution $X\in\C^{n\times n}$ of $\exp X=Q^*D$ and that $\tr{\sym\Log Q^* D}$ is real.

For any $Q\in \U(n)$ the Hermitian positive definite matrices $\exp(\sym\Log Q^* D)$ and $\exp(-\sym\Log Q^* D)$ can be simultanuously unitarily diagonalized with positive eigenvalues, i.e., for some $Q_1\in\U(n)$
\begin{align}
  \label{eigenvalue_representation}
    Q_1^* \exp(\sym\Log Q^* D) Q_1&=\exp(Q_1^* (\sym\Log Q^* D) Q_1)=\diag(x_1,\ldots,x_n)   \, ,\notag\\
Q_1^*\exp(-\sym\Log Q^* D) Q_1&=\exp(-Q_1^* (\sym\Log Q^* D) Q_1)  \notag\\
  &=\left(\exp(Q_1^* (\sym\Log Q^* D) Q_1\right)^{-1}
   =  \diag(\frac{1}{x_1},\ldots, \frac{1}{x_n})\, ,
   \end{align}
since $X\mapsto \exp X$ is an isotropic function. We arrange the positive real eigenvalues in decreasing order 
 $x_1\ge x_2\ge\ldots \ge x_n>0$. For any unitarily invariant norm it follows therefore from \eqref{exp_estimate_1}, \eqref{exp_estimate_2} and \eqref{det_equality} together with \eqref{eigenvalue_representation} that
\begin{align}
\label{unitary_matrix_estimate}
&\norm{\diag(x_1,\ldots,x_n)}^2\!
=\!\norm{Q_1^* \exp(\sym\Log Q^* D) Q_1}^2\!=\!\norm{\exp(\sym\Log Q^* D) }^2 \!\ge \!\!\norm{D}^2  \\   
&\norm{  \diag(\frac{1}{x_1},\ldots, \frac{1}{x_n})}^2 =\norm{Q_1^* \exp(-\sym\Log Q^* D) Q_1}^2=\norm{\exp(-\sym\Log Q^* D) }^2 \ge \norm{D^{-1}}^2 \notag\\
& \Det\diag(x_1,\ldots,x_n)=\det{Q_1^*\exp(\sym\Log Q^* D) Q_1}=\det{\exp(\sym\Log Q^* D)}=\det D\, . \notag
\end{align}
Below we combine these inequalities and the ``sum of squared logarithms inequality'' to give a proof of Theorem \ref{thm:main}.
\subsection{Frobenius matrix norm for $n=2,3$}
Now consider the Frobenius matrix norm for dimension $n=3$. The three conditions in \eqref{unitary_matrix_estimate} can be expressed as
\begin{align}
\label{algebraic_estimate}
     x_1^2+x_2^2+x_3^2&\ge d_1^2+d_2^2+d_3^2  \notag \\
     \frac{1}{x_1^2}+\frac{1}{x_2^2}+\frac{1}{x_3^2}&\ge \frac{1}{d_1^2}+\frac{1}{d_2^2}+\frac{1}{d_3^2}\\
     x_1 \, x_2\, x_3&= d_1\, d_2\, d_3. \notag
\end{align}
By a new result: the ``sum of squared logarithms inequality'' \cite{Neff_log_inequality13}, conditions \eqref{algebraic_estimate} imply 
\begin{align}
 \label{squared_log_inequality}
   (\log x_1)^2+(\log x_2)^2+(\log x_3)^2\ge  (\log d_1)^2+(\log d_2)^2+(\log d_3)^2\, ,
\end{align}
with equality if and only if $(x_1,x_2,x_3)=(d_1,d_2,d_3)$. 
This is true, despite the map $t\mapsto (\log t)^2$ being non-convex. 
Similarly, for the two-dimensional case with a much simpler proof \cite{Neff_log_inequality13} 
\begin{align}
\label{algebraic_estimate_2d}
\left.\begin{array}{ll}     x_1^2+x_2^2&\ge d_1^2+d_2^2 \\
\displaystyle     \frac{1}{x_1^2}+\frac{1}{x_2^2}&\displaystyle \ge \frac{1}{d_1^2}+\frac{1}{d_2^2}  \\
     x_1 \, x_2 &= d_1\, d_2  
\end{array}\right\}
\quad
\Rightarrow \quad (\log x_1)^2+(\log x_2)^2\ge  (\log d_1)^2+(\log d_2)^2\, .
\end{align}
Since on the one hand \eqref{log_on_psym_1} and \eqref{log_on_psym_2} imply
\begin{align}
\label{X_estimate}
  (\log x_1)^2+(\log x_2)^2+(\log x_3)^2&=\norm{\logp \diag(x_1,x_2,x_3)}_F^2 \notag\\
 &=\norm{\logp (Q_1^* \exp(\sym\Log Q^* D) Q_1)}_F^2   \\
 &=\norm{Q_1^*   \logp \exp(\sym\Log Q^* D) Q_1}_F^2  \notag\\
 &=\norm{ \logp \exp(\sym\Log Q^* D)}_F^2=\norm{\sym\Log Q^* D}_F^2 \notag
 \end{align}
and clearly
\begin{align}
\label{D_equality}
  (\log d_1)^2+(\log d_2)^2+(\log d_3)^2=\norm{\logp D}_F^2\, ,
\end{align}
we may combine \eqref{X_estimate} and \eqref{D_equality} with the sum of squared logarithms inequality \eqref{squared_log_inequality} to obtain
\begin{align}\label{eq:symD}
   \norm{\sym\Log Q^* D}_F^2\ge \norm{\logp D}_F^2
   \end{align}
for any $Q\in\U(3)$. Since on the other hand we have the trivial upper bound (choose $Q=\id$)
\begin{align}
\min_{Q\in \U(3)}\norm{\sym \Log (Q^* D) }_F^2\le \norm{\logp D}_F^2\, ,
\end{align}
this shows that
\begin{align}
\min_{Q\in \U(3)}\norm{\sym \Log (Q^* D) }_F^2=\norm{\logp D}_F^2\, .
\end{align}
The minimum is realized for $Q=\id$, which corresponds to the polar factor $U_p$ in the original formulation. Noting that
\begin{align}\label{eq:skewherm}
\norm{\Log(Q^* D)}_F^2= \norm{\sym \Log (Q^* D) }_F^2+ \norm{\skew \Log (Q^* D) }_F^2\ge \norm{\sym \Log (Q^* D) }_F^2  
\end{align}
 by the orthogonality of the Hermitian and skew-Hermitian parts in the trace scalar product, we also obtain 
\begin{align}\label{eq:frob3}
  \min_{Q\in \U(3)}\norm{\Log (Q^* D) }_F^2\ge \min_{Q\in \U(3)}\norm{\sym \Log (Q^* D) }_F^2=\norm{\logp D}_F^2\, .
  \end{align}
Since all the terms in \eqref{eq:frob3} are equal when $Q=I$ and the principal logarithm is taken, we obtain
\begin{align}\label{eq:frob4}
  \min_{Q\in \U(3)}\norm{\Log (Q^* D) }_F^2=\norm{\logp D}_F^2\, .
  \end{align}
Hence, combining again we obtain for all $\mu>0$ and all $\mu_c\ge 0$
\begin{align}
\min_{Q\in \U(3)}  \mu\, \norm{\sym \Log (Q^* D) }_F^2+\mu_c\, \norm{\skew \Log (Q^* D) }_F^2 =\mu\, \norm{\logp D}_F^2\, .
\end{align}

Observe that although we allowed $\Log$ to be any matrix logarithm, the one that gives the smallest $\norm{\Log (Q^* D) }_F$ and $\norm{\sym \Log (Q^* D) }_F$ is in both cases the principal logarithm, regardless of $Z$. 

\subsection{Spectral matrix norm for arbitrary $n\in\N$}
For the spectral  norm, the conditions \eqref{unitary_matrix_estimate} can be expressed as 
\begin{align}
  \label{spectral_conditions}
     x_1^2\ge d_1^2\, , & \quad
     \frac{1}{x_n^2}\ge \frac{1}{d_n^2}\, ,\notag\\
     x_1 \, x_2\, x_3\, \ldots x_n&= d_1\, d_2\, d_3\, \ldots d_n.  
\end{align}
This yields the ordering
\begin{align}
    0<x_n \le  d_n  \le d_1\le  x_1\, .
\end{align}
It is easy to see that this implies (even without the determinant condition \eqref{spectral_conditions})
\begin{align}
     \max \{ \abs{\log x_n}\, , \abs{\log d_n}\, , \abs{\log d_1}\, , \abs{\log x_1}\}= \max \{ \abs{\log x_n}\, , \abs{\log x_1}\}    \, ,
\end{align}
which shows 
\begin{align}
     \max \{ \abs{\log d_n}\, , \abs{\log d_1}\} \le \max \{ \abs{\log x_n}\, , \abs{\log x_1}\} \,.   
\end{align}
Therefore, cf. \eqref{X_estimate},
\begin{align}
\norm{\sym\Log Q^* D}_2^2
=&\norm{\logp \diag(x_1, \ldots, x_n)}_{2}^2\notag\\
=&\norm{\diag(\log x_1, \ldots, \log x_n)}_{2}^2\notag\\
=&\max_{i=1,2,3,\, \ldots , n}\{ \abs{\log x_1}\, , \abs{\log x_2}\, , \ldots\, , \abs{\log x_n}\}^2 \notag\\
=&\max_{i=1,2,3,\, \ldots , n}\{ (\log x_1)^2\, , (\log x_2)^2\, , \ldots\, , (\log x_n)^2\} \notag\\
=&\max \{ (\log x_1)^2\, , (\log x_n)^2\} \notag\\
&\ge \max \{ (\log d_1)^2\, , (\log d_n)^2\}\notag\\
&=\max_{i=1,2,3,\ldots n}\{ (\log d_1)^2\, , (\log d_2)^2\, , \ldots\, , (\log d_n)^2\}\\
&=\max_{i=1,2,3,\, \ldots , n}\{ \abs{\log d_1}\, , \abs{\log d_2}\, , \ldots\, , \abs{\log d_n}\}^2 \notag\\
&=\norm{\diag(\log d_1, \ldots, \log d_n)}_{2}^2\notag\\
&= \norm{\log\diag(d_1, \ldots, d_n)}_{2}^2 = \norm{\logp D}_{2}^2\, ,\notag
\end{align}
from which we obtain, as in the case of the Frobenius norm, due to unitary invariance, 
\begin{align}
\min_{Q\in \U(n)}\norm{\sym \Log (Q^* D) }_{2}^2=\norm{\logp D}_{2}^2\, .
\end{align}
For complex numbers we have the bound $\abs{z}\ge \abs{\Re z}$. A matrix analogue is that the spectral norm of some matrix $X\in\C^{n\times n}$ bounds the spectral norm of the Hermitian part $\sym X$, see \cite[p.355]{Bernstein2009} and \cite[p.151]{Horn91}, i.e.
$
    \norm{X}_2^2\ge \norm{\sym X}_2^2
$. In fact, this inequality holds for all unitarily invariant norms \cite[p.454]{Horn85}:
\begin{align}
\label{sym_estimate}
 \forall\, X\in \C^{n\times n}:\quad  \norm{X}^2\ge \norm{\sym X}^2\, .
  \end{align}
Therefore we conclude that for the spectral norm, in any dimension we have
\begin{align}
  \min_{Q\in \U(n)}\norm{\Log (Q^* D) }_2^2\ge \min_{Q\in \U(n)}\norm{\sym \Log (Q^* D) }_2^2=\norm{\logp D}_2^2\, ,
  \end{align}
with equality holding for $Q=U_p$. 

\section{The real Frobenius case on $\SO(3)$}\label{sec:real}
In this section we consider $\A\in\GL^+(3,\R)$, which implies that $\A=U_p \, H$ admits the polar decomposition with $U_p\in\SO(3)$ and 
an eigenvalue decomposition $H=VDV^T$ for  $V\in \SO(3)$. 
We observe that
\begin{align}
\label{real_estimate}
    \min_{Q\in \SO(n)}\norm{\sym \Log (Q^T D) }_F^2\ge  \min_{Q\in \U(n)}\norm{\sym \Log (Q^* D) }_F^2   \, .
\end{align}
Therefore, for all $\mu>0, \, \mu_c\ge 0$ we have, using inequality \eqref{real_estimate}
  \begin{align}
    \min_{Q\in\SO(3)} \mu\, &\| \sym \Log(Q^T \A)\|_F^2+\mu_c\, \|\skew \Log(Q^T \A)\|_F^2 \notag\\
        &\ge  \min_{Q\in\SO(3)} \mu\, \| \sym \Log(Q^T \A)\|_F^2 \\
   & =\mu \norm{\sym \log(U_p^T \A)}_F^2  
   =\mu \norm{\sym \log(U_p^T \A)}_F^2+  \mu_c \underbrace{ \norm{\skew \log(U_p^T \A)}_F^2}_{\text{$=0$}} \, ,  \notag
     \end{align} 
and it follows that the solution to the minimization problem \eqref{general_problem} for $\A\in \GL^+(n,\R)$ and $n=2,3$ is also obtained by the orthogonal polar factor (a similar argument holds for $n=2$).

Denoting by $\dev_n X=X-\frac{1}{n}\tr{X}\id$ the orthogonal projection of $X\in\R^{n\times n}$ onto trace free matrices in the trace scalar product, 
we obtain a further result of interest in its own right (in which we really need $Q\in\SO(3)$), namely
\begin{align}
\label{hencky_isochoric}
 \min_{Q\in \SO(3)}\norm{\dev_3\Log (Q^T D) }_F^2& =\norm{\dev_3\log D }_F^2\, ,\notag\\
  \min_{Q\in \SO(3)}\norm{\dev_3\sym \Log (Q^T D) }_F^2 & =\norm{\dev_3\log D }_F^2\, .
\end{align}
As was in the previous section, it suffices to show the second equality.
This is true since by using \eqref{log_det_computation} for $Q\in\SO(3)$ we have
\begin{align}
  \min_{Q\in \SO(3)}\norm{\dev_3\sym \Log (Q^T D) }_F^2
    &=\min_{Q\in \SO(3)}\left(\norm{\sym \Log (Q^T D) }_F^2-\frac{1}{3}\tr{\Log Q^T D}^2\right)\notag\\
    &=\min_{Q\in \SO(3)}\left(\norm{\sym \Log (Q^T D) }_F^2-\frac{1}{3}(\log \det{Q^T D})^2\right) \notag\\    
    &=\min_{Q\in \SO(3)}\norm{\sym \Log (Q^T D) }_F^2-\frac{1}{3}(\log \det{D})^2 \notag\\    
    &=\min_{Q\in \SO(3)}\norm{\sym \Log (Q^T D) }_F^2-\frac{1}{3}\tr{\log D}^2\notag\\
  &\ge \min_{Q\in \U(n)}\norm{\sym \Log (Q^* D) }_F^2-\frac{1}{3}\tr{\log D}^2\notag\\
    &=\norm{\sym \log D }_F^2-\frac{1}{3}\tr{\log D}^2\\   
        &=\norm{\sym \log D }_F^2-\frac{1}{3}\tr{\sym\log D}^2\notag\\       
    &=\norm{\dev_3\sym \log D }_F^2=\norm{\dev_3\log D }_F^2 \, .\notag
\end{align}

\section{Uniqueness}\label{sec:unique}
We have seen that the polar factor $U_p$ minimizes both $\|\Log(Q^* \A)\|^2$ and $\|\sym(\Log(Q^* \A))\|^2$, but what about its uniqueness? Is there any other unitary matrix that also attains the minimum? 
We address these questions below. 

\subsection{Uniqueness of $U_p$ as the minimizer of $\|\Log(Q^* \A)\|^2$}\label{sec:unique1}
Note that the unitary polar factor $U_p$ itself is not unique when $\A$ does not have full column rank~\cite[Thm.~8.1]{Higham2008}. However in our setting we do not consider this case because $\Log(U\A)$ is defined only if $U\A$ is nonsingular. 

We show below that $U_p$ is the unique minimizer of $\|\Log(Q^* \A)\|^2$ for the Frobenius norm, while for the spectral norm there can be many $Q\in \U(n)$ for which $\|\log(Q^* \A)\|^2=\|\log(U_p^* \A)\|^2$. 
\paragraph{Frobenius norm for $n \leq 3$.}
We focus on $n=3$ as the case $n=2$ is analogous and simpler. 
By the fact that $Q=U_p$ satisfies equality in \eqref{eq:frob3}, any minimizer $Q$ of 
$\norm{\Log (Q^* D) }_F$ must satisfy 
\begin{align}\label{eq:need}
\norm{\Log (Q^* D) }_F= \norm{\sym \Log (Q^* D) }_F=\norm{\log D}_F. 
\end{align}
Note that by \eqref{eq:skewherm} the first equality of \eqref{eq:need} holds only if $\Log(Q^* D)$ is Hermitian.

We now examine the condition that satisfies the latter equality of \eqref{eq:need}. 
Since 
$\Log (Q^* D)$ is Hermitian the matrix $\exp(\Log (Q^* D))$ is positive definite, 
so we can write $\exp(\Log (Q^* D)) = Q_1^*\diag(x_1,x_2,x_3)Q_1$ 
 for some unitary $Q_1$ and $x_1,x_2,x_3>0$. Therefore 
 \begin{align}\label{eq:qd}
\log (Q^* D) = Q_1^*\diag(\log x_1,\log x_2,\log x_3)Q_1.    
 \end{align}
Hence for $\norm{\sym \Log (Q^* D) }_F=\norm{\log D}_F$ to hold we need
\[    (\log x_1)^2+(\log x_2)^2+(\log x_3)^2=  (\log d_1)^2+(\log d_2)^2+(\log d_3)^2,\]
which is precisely the case where equality holds in the sum of squared logarithms inequality \eqref{squared_log_inequality}. 
As discussed above, equality holds in \eqref{squared_log_inequality} if and only if 
 $(x_1,x_2,x_3)=(d_1,d_2,d_3)$. Hence by \eqref{eq:qd} we have $\log (Q^* D) = Q_1^*\diag(\log x_1,\log x_2,\log x_3)Q_1 = Q_1^*\log(D)Q_1$, so 
taking the exponential on both sides yields
\begin{align}\label{eq:qd2}
Q^* D = Q_1^*DQ_1.  
\end{align}
Hence $Q_1Q^* DQ_1^* = D$. Since $Q_1Q^*$ and $Q_1^*$ are both unitary matrices this is a singular value decomposition of $D$. 
Suppose $d_1>d_2>d_3$. 
Then since the singular vectors of distinct singular values are unique up to multiplication by $e^{i\vartheta}$, it follows that 
$Q_1Q^* = Q_1 = \diag(e^{i\vartheta_1},e^{i\vartheta_2},e^{i\vartheta_3})$ for $\vartheta_i\in\mathbb{R}$, 
so $Q= I$. If some of the $d_i$ are equal, for example if $d_1=d_2>d_3$, then we have 
$Q_1 = \diag(Q_{1,1}, e^{i\vartheta_3})$ where $Q_{1,1}$ is a $2\times 2$ arbitrary unitary matrix, but we still have $Q = I$. 
If $d_1=d_2=d_3$, then $Q_1$ can be any unitary matrix but again $Q = I$. 
Overall,  for \eqref{eq:qd2} to hold we always need $Q = I$, which corresponds to the unitary polar factor $U_p$ in the original formulation. 
Thus $Q = U_p$ is the unique minimizer of $\norm{\Log (Q^* D) }_F$ with minimum $\norm{\log (U_p^* D) }_F$. 

\paragraph{Spectral norm}
For the spectral norm there can be many unitary matrices $Q$ that attain $\|\Log(Q^*\A)\|_2^2=\|\log(U_p^*\A)\|_2^2$. For example, consider 
$\A=\big[\begin{smallmatrix}  e&0\cr  0&1\end{smallmatrix}\big]$. The unitary polar factor is $U_p=\id$. 
Defining $U_1=\big[\begin{smallmatrix}  1&0\cr  0&e^{i\vartheta}\end{smallmatrix}\big]$
we have $\|\log(U_1\A)\|_2=\|\big[\begin{smallmatrix}  1&0\cr  0&\vartheta\end{smallmatrix}\big]\|_2=
1$ for any $\vartheta\in[-1,1]$. 

Now we discuss the general form of the minimizer $Q$. Let $\A=U\Sigma V^*$ be the SVD with 
$\Sigma = \diag(\sigma_1,\sigma_2,\ldots, \sigma_n)$. 
Recall that $\|\log (U_p^*\A)\|_2 = \max(|\log\sigma_1(\A)|,|\log\sigma_n(\A)|)$. 

Suppose that $\|\log (U_p^*\A)\|_2 = |\log\sigma_1(\A)|\geq |\log\sigma_n(\A)|$. Then 
for any $Q = U\diag(1,Q_{22})V^*$ we have 
\[\log Q^*\A =  \log V \diag(1,Q_{22} )\Sigma V^*, \]
so we have $\|\log Q^*\A\|_2 = |\log\sigma_1(\A)|= \|\log U_p^*\A\|_2$ for any  $Q_{22}\in \U(n-1)$ such that 
$\|\log Q_{22}\diag(\sigma_2,\sigma_3,\ldots, \sigma_n)\|_2\leq \|\log U_p^*\A\|_2 $. Note that such $Q_{22}$ always includes $I_{n-1}$, but 
may not include the entire set of $(n-1)\times (n-1)$ unitary matrices as evident from the above simple example. 

Similarly, if $\|\log U_p^*\A\|_2 = |\log\sigma_n(\A)|\geq |\log\sigma_1(\A)|$, then 
we have $\|\log Q^*\A\|_2 = \|\log U_p^*\A\|_2$ for 
 $Q = U\diag(Q_{22},1)V^*$ where  $Q_{22}$ can be any $(n-1)\times (n-1)$ unitary matrix satisfying 
$\|\log Q_{22}\diag(\sigma_1,\sigma_2,\ldots, \sigma_{n-1})\|_2\leq \|\log U_p^*\A\|_2 $.

\subsection{Non-uniqueness of $U_p$ as the minimizer of $\|\sym(\Log(Q^* \A))\|^2$}
The fact that $U_p$ is not the unique minimizer of $\|\sym(\Log(Q^* \A))\|^2$ can be seen by the simple example $\A=\id$. Then 
$\Log Q^*$ is a skew-Hermitian matrix, so $\sym(\Log(Q^* \A))=0$ for \emph{any} unitary $Q$.

In general, every $Q$ of the following form gives the same value of $\|\mbox{sym}(\Log(Q^* \A))\|^2$. 
Let $\A=U\Sigma V^*$ be the SVD with $\Sigma=\mbox{diag}(\sigma_1I_{n_1},\sigma_2I_{n_2},\ldots, \sigma_kI_{n_k})$ where $n_1+n_2+\cdots+n_k$ ($k$ if $\A$ has pairwise distinct singular values). Then it can be seen that any unitary $Q$ of the form 
\begin{align}\label{eq:nonunique}
Q^*=U\mbox{diag}(Q_{n_1},Q_{n_2},\ldots, Q_{n_k})V^*, 
\end{align}
where $Q_{n_i}$ is any $n_i\times n_i$ unitary matrix, 
yields $\|\sym(\Log(Q^* \A))\|^2=\|\sym(\log(U_p^*\A))\|^2$. Note that this holds for any unitarily invariant norm. 

The above argument naturally leads to the question of whether $U_p$ is unique up to $Q_{n_i}$ in \eqref{eq:nonunique}. In particular, when the singular values of $Z$ are distinct, is $U_p$ determined up to scalar rotations $Q_{n_i} = e^{i\vartheta_{n_i}}$? 

For the spectral norm an argument similar to that above shows there can be many $Q$ for which $\|\sym(\Log(Q^* \A))\|^2=\|\sym(\log(U_p^*\A))\|^2$. 

For the Frobenius norm, the answer is yes. To verify this, observe in \eqref{eq:symD} that $\norm{\sym \Log (Q^* D) }_F=\norm{\log D}_F$ implies  $(x_1,x_2,x_3)=(d_1,d_2,d_3)$ and hence 
$\Log (Q^* D)
= Q_1^*\diag(\log d_1,\log d_2,\log d_3)Q_1+S$, where $S$ is a skew-Hermitian matrix. Hence 
\begin{equation}  \label{eq:swhat}
\exp(Q_1^*\diag(\log d_1,\log d_2,\log d_3)Q_1+S) = Q^*D,   
\end{equation}
and by \eqref{bhatia} we have
\begin{align*}
\|Q^*D\|_F& = \|\exp(Q_1^*\diag(\log d_1,\log d_2,\log d_3)Q_1+S)\|_F\\
&\leq \|\exp(Q_1^*\diag(\log d_1,\log d_2,\log d_3)Q_1)\|_F=\|Q^*D\|_F. 
\end{align*}
Since equality in \eqref{bhatia} holds for the Frobenius norm if and only if $X$ is normal  (which can be seen from the proof of~\cite[Thm.~IX.3.1]{Bhatia97}), 
for the last inequality to be an equality, 
$Q_1^*\diag(\log d_1,\log d_2,\log d_3)Q_1+S$ must be a normal matrix. 
Since $Q_1^*\diag(\log d_1,\log d_2,\log d_3)Q_1$ is Hermitian and $S$ is skew-Hermitian, this means $Q_1^*\diag(\log d_1,\log d_2,\log d_3)Q_1+S=Q_1^*\diag(is_1+\log d_1,is_2+\log d_2,is_3+\log d_3)Q_1$ for $s_i\in\R$. 
Together with \eqref{eq:swhat} we conclude that 
\[Q^*D = Q_1^*\diag(d_1e^{is_1},d_2e^{is_2},d_3e^{is_3})Q_1.\]
By an argument similar to that following \eqref{eq:qd2} we obtain $Q=\diag(e^{-is_1},e^{-is_2},e^{-is_3})$.

\section{Conclusion and outlook}\label{sec:conclude}
The result in the Frobenius matrix norm cases for $n=2,3$ hinges crucially on the use of the new sum of squared logarithms inequality \eqref{squared_log_inequality}. This inequality seems to be true in any dimensions with appropriate additional conditions \cite{Neff_log_inequality13}. However, we do not have a proof yet.

Nevertheless, numerical experiments suggest that the optimality of the polar factor $U_p$ in both
\begin{align}
 \min_{Q\in\U(n)} \| \Log(Q^* \A)\|^2\, ,\quad  \min_{Q\in\U(n)} \| \sym \Log(Q^* \A)\|^2
  \end{align}
is true for any unitarily invariant norm, over $\R$ and $\C$ and in any dimension. This would imply that for all $\mu,\mu_c\ge 0$ and for any unitarily invariant norm
\begin{align*}
   \min_{Q\in\U(n)} \mu\, \| \sym \Log(Q^* \A)\|^2+\mu_c\,\norm{\skew \Log(Q^* \A)}^2=\mu\, \norm{\log(U_p^* \A) }^2 =\mu\,\norm{\log\sqrt{\A^* \A}}^2  \, .
   \end{align*}
 We also conjecture that $Q=U_p$ is the unique unitary matrix that minimizes $\|\Log(Q^* \A)\|^2$ for every unitarily invariant norm.

In a forthcoming contribution \cite{Neff_Osterbrink_hencky13} we will use our new characterization of the orthogonal factor in the polar decomposition to calculate the geodesic distance of the isochoric part of the deformation gradient $ \frac{F}{\det{F}^{\frac{1}{3}}}    \in\SL(3)$ to $\SO(3)$ in the canonical left-invariant Riemannian metric on $\SL(3)$,
namely based on \eqref{hencky_isochoric}
\begin{align}
    \dist_{{\rm geod},\SL(3)}^2(\frac{F}{\det{F}^{\frac{1}{3}}},\SO(3) )=\norm{\dev_3\log\sqrt{F^T F}}_F^2=\min_{Q\in\SO(3)}\norm{\dev_3\sym\Log Q^T F}_F^2\, . \notag
\end{align}
Thereby, we provide a rigorous geometric justification for the preferred use of the Hencky-strain measure $\norm{\log\sqrt{F^TF}}_F^2$ in nonlinear elasticity and plasticity theory, see \cite{Hencky28,Xiao2005} and the references therein.

\subsection*{Acknowledgment}
The authors are grateful to N. J. Higham (Manchester) for establishing their contact and for valuable remarks. P. Neff acknowledges a helpful computation of J. Lankeit (Essen), which led to the crucial formulation of estimate \eqref{exp_estimate_2}. We thank the referees for their remarks and suggestions.

\section{Appendix}
\subsection{Connections between $\C^\times$ and $\CO(2)$}
The following connections between $\C^\times$ and $\CO(2)=\R^+\cdot \SO(2)$ are clear:
\begin{align}
    \abs{z}^2=a^2+b^2\quad & =\quad  \norm{\A}_{\CO}^2=\frac{1}{2}\norm{\A}_F^2 \, ,\notag\\
    \overline{z}=a-i b=a+i (-b) \quad &\Leftrightarrow\quad  \A^T=\begin{bmatrix}
     a & -b\\
     b & a
 \end{bmatrix} \, ,\notag\\   
     & \quad \A\,\A^T=\A^T \A
     =\begin{bmatrix}
          a^2+b^2 & 0 \\
          0 & a^2+b^2
      \end{bmatrix}   \, ,     \\
       \overline{z}\,z=\abs{z}^2 \quad &=\quad  \norm{\A}_{\CO}^2=\frac{1}{2}\tr{\A^T \A}=\frac{1}{2}\norm{\A}_F^2\, ,\notag\\    
    z\cdot w=w\cdot z \quad &\Leftrightarrow\quad  \A\cdot W=W\cdot \A\, ,\notag
\end{align}
\begin{align}    
     \Re(z)=\frac{z+\overline{z}}{2}=a \quad &\Leftrightarrow\quad  \sym(\A)=\frac{1}{2}(\A+\A^T)=\begin{bmatrix}
     a & 0\\
     0 & a
 \end{bmatrix} \, ,\notag\\
    \Im(z)=\frac{z-\overline{z}}{2}=b \quad &\Leftrightarrow\quad  \skew(\A)=\frac{1}{2}(\A-\A^T)=\begin{bmatrix}
     0 & b\\
     -b & 0
 \end{bmatrix} \, ,   \\
    \abs{\Re(z)}^2=\abs{a}^2   \quad &=\quad  \norm{\sym\A}_{\CO}^2=\frac{1}{2}\norm{\sym \A}_F^2 \, ,\notag\\
    \abs{\Im(z)}^2=\abs{b}^2   \quad &=\quad  \norm{\skew\A}_{\CO}^2=\frac{1}{2}\norm{\skew \A}_F^2 \, ,\notag\\
\abs{z}^2=\Re(z)^2+\Im(z)^2   \quad & =\quad  \det{\A}\, ,\notag
\end{align}
\begin{align}
e^{i\vartheta}=\cos\vartheta+i \sin\vartheta\quad &\Leftrightarrow\quad  
Q(\vartheta)=
\begin{bmatrix} 
\cos\vartheta &\sin\vartheta\\
-\sin\vartheta & \cos\vartheta 
\end{bmatrix}\in\SO(2)
\, ,\quad \vartheta\in (-\pi,\pi]\, ,\notag\\
e^{-i\vartheta} \,z =  (e^{i\vartheta})^{-1} \, z \quad &\Leftrightarrow\quad  Q^T \A \, ,\notag\\
 z=e^{i\arg(z)}\, \abs{z} \quad &\Leftrightarrow\quad  \A=U_p\, H\, ,\quad \text{{\bf polar form} versus {\bf polar decomposition}} \notag\\
 \abs{z} \quad & \Leftrightarrow\quad \; H=\sqrt{\A^T\A}=\sqrt{\det{\A}}\id_2\, ,\notag\\
 & \quad \quad \;U_p =\A H^{-1}=\frac{1}{\sqrt{a^2+b^2}}
 \begin{bmatrix}
   a & b\\
   -b & a
\end{bmatrix}\in\SO(2) \, ,  \\
\abs{e^z}= \abs{e^{a+ib}}=\abs{e^{\Re(z)}}\quad & \Leftrightarrow\quad \norm{\exp(\A)}_F=\norm{\exp(a\id_2+\skew(\A))}_F\notag\\
&\quad =\norm{\exp(a\id_2)\,\exp(\skew(\A))}_F=\norm{\exp(a\id_2)}_F=\norm{\exp(\sym(\A))}_F\, .\notag
\end{align}

\subsection{Optimality properties of the polar form}
The polar decomposition is the matrix analog of the polar form of a complex number
\begin{align}
      z=e^{i\arg(z)} \abs{z}\, .
\end{align}
The argument $\arg(z)$ determines the unitary part $e^{i\arg(z)}$ while the positive definite Hermitian matrix is $\abs{z}$. 
The argument $\arg(z)$ in the polar form is optimal in the sense that
\begin{align}
 \min_{\vartheta\in (-\pi,\pi]} \mu\, \abs{e^{-i\vartheta}z-1}^2=\min_{\vartheta\in (-\pi,\pi]} \mu\, \abs{z-e^{i\vartheta}}^2=\mu\, \abs{\abs{z}-1}^2\,, \quad \vartheta=\arg(z) \, .
 \end{align}
 However, considering only the real (Hermitian) part
 \begin{align}
      \min_{\vartheta\in (-\pi,\pi]} \mu\, \abs{\Re (e^{-i\vartheta}z-1) }^2
 =\begin{cases}
              \mu\, \abs{\abs{z}-1}^2 & \quad \abs{z}\le 1\, ,\quad \vartheta=\arg(z)\\
              0 & \quad \abs{z}>1\, , \quad \vartheta:\; \cos(\arg(z)-\vartheta)=\frac{1}{\abs{z}}\, ,
  \end{cases} 
  \end{align}
  shows that optimality of $\vartheta=\arg(z)$ ceases to be true for $\abs{z}>1$ and the optimal $\vartheta$ is not unique. This is the \emph{nonclassical solution} alluded to in \eqref{symmetry_break}. In fact we have optimality of the polar factor for the \emph{Euclidean weighted family} only for $\mu_c\ge \mu$:
  \begin{align}
   \min_{\vartheta\in (-\pi,\pi]} \mu\, \abs{\Re(e^{-i\vartheta}z-1)}^2+\mu_c\, \abs{\Im(e^{-i\vartheta}z-1)}^2
   = \mu\, \abs{\abs{z}-1}^2\, ,\quad \vartheta=\arg(z)\, ,
  \end{align}
  while for $\mu\ge \mu_c\ge 0$ there always exists a $z\in\C$ such that 
   \begin{align}
   \min_{\vartheta\in (-\pi,\pi]} \mu\, \abs{\Re(e^{-i\vartheta}z-1)}^2+\mu_c\, \abs{\Im(e^{-i\vartheta}z-1)}^2
     <\mu\, \abs{\abs{z}-1}^2\, .
  \end{align} 
  In pronounced contrast, for the \emph{logarithmic weighted family} the polar factor is optimal for all choices of weighting factors $\mu,\mu_c\ge 0$:
  \begin{align}
     \min_{\vartheta\in (-\pi,\pi]} \mu\, \abs{\Re\Log_\C(e^{-i\vartheta}z)}^2+\mu_c\, \abs{\Im\Log_\C(e^{-i\vartheta}z)}^2=
 \begin{cases}    
     \mu\, \abs{\log|z|}^2   &\quad \mu_c>0: \;\vartheta=\arg(z)\\  
      \mu\, \abs{\log|z|}^2  &\quad \mu_c=0: \; \vartheta \;\text{arbitrary}\, .
      \end{cases}
\end{align}
Thus we may say that the more fundamental characterization of the polar factor as minimizer is given by the property with respect to the logarithmic weighted family.

\subsection{The three-parameter case $\SL(2)$ by hand}

\subsubsection{Closed form exponential on $\SL(2)$ and closed form principal logarithm}
The exponential on $\sL(2)$ can be given in closed form, see \cite[p.78]{Tarantola06} and \cite{Bernstein93}.
Here, the two-dimensional Caley-Hamilton theorem is useful: 
$X^2-\tr{X}X+\det{X}\id_2=0$. Thus, for $X$ with $\tr{X}=0$, it holds $X^2=-\det{X}\id_2$ and $\tr{X^2}=-2\,\det{X}$. Moreover, every higher exponent $X^k$ can be expressed in $\id$ and $X$ which shows that $\exp(X)=\alpha(X)\id_2+\beta(X) X$. 
Tarantola \cite{Tarantola06} defines the "near zero subset" of $\sL(2)$
\begin{align}
   \sL(2)_0:= \{X\in\sL(2)\, |\, \Im\left(\sqrt{\frac{1}{2}\tr{X^2}}\right) < \pi\, \}
\end{align}
and the "near identity subset" of $\SL(2)$
  \begin{align}
   \SL(2)_{\id}:= &\SL(2)\setminus\{\text{both eigenvalues are real and negative}\, \}\, .
\end{align}
On this set the principal matrix logarithm is real. Complex eigenvalues appear always in conjugated pairs, therefore the eigenvalues are either real or complex in the twodimensional case. Then it holds \cite[p.149]{Helgason01}
\begin{align}
\label{exp_formula_2d}
  \exp(X)=\begin{cases}
        \cosh(\sqrt{-\det{X}})\, \id_2+ \frac{\sinh(\sqrt{-\det{X}})}{\sqrt{-\det{X}}}\, X     & \det{X}<0\\
         \cos(\sqrt{\det{X}})\, \id_2+  \frac{\sin(\sqrt{\det{X}})}{\sqrt{\det{X}}}\, X    & \det{X}>0\\
         \id_2+X                & \det{X}=0\, .
\end{cases}
\end{align}
Therefore
\begin{align}
\label{exp_sl2}
   \exp:&\; \sL(2)_0\mapsto\SL(2)_\id\, ,\, 
   \exp(X)=\cosh(s)\id_2+\frac{\sinh(s)}{s}\, X\, ,\quad s:=\sqrt{\frac{1}{2}\tr{X^2}}
    =\sqrt{-\det{X}}\, .
    \end{align}
 Since the argument $s=\sqrt{-\det{X}}$ is complex valued  for $\det{X}>0$ we note that $\cosh(i y)=\cos(y)$.

The one-parameter $\SO(2)$ case is included in the former formula for the exponential. The previous formula \eqref{exp_sl2} can be specialized to $\so(2,\R)$. Then
\begin{align}
\exp\begin{bmatrix}
         0 & \alpha\\
         -\alpha & 0
\end{bmatrix}&=\cosh(\sqrt{-\alpha^2})\id_2
+\frac{\sinh(\sqrt{-\alpha^2})}{\sqrt{-\alpha^2}}
\begin{bmatrix}
      0 & \alpha\\
      -\alpha & 0
\end{bmatrix}
  =\cosh(i \abs{\alpha})\id_2
+\frac{\sinh(i \abs{\alpha})}{i\abs{\alpha}}
\begin{bmatrix}
      0 & \alpha\\
      -\alpha & 0
\end{bmatrix}\notag\\
  &=\cos(\abs{\alpha})\id_2
+\frac{i\sin(\abs{\alpha})}{i\abs{\alpha}}
\begin{bmatrix}
      0 & \alpha\\
      -\alpha & 0
\end{bmatrix}
=
\begin{bmatrix}
      \cos\alpha & \sin\alpha\\
      -\sin\alpha & \cos\alpha
\end{bmatrix}\, .
\end{align}


We need to observe that the exponential function is not surjective onto $\SL(2)$ since not every matrix $\A\in\SL(2)$ can be written as $\A=\exp(X)$ for $X\in\sL(2)$. This is the case because $\forall\;X\in\sL(2):\quad  \tr{\exp(X)}\ge -2$, see $\eqref{exp_formula_2d}_2$. Thus, any matrix $\A\in\SL(2)$ with $\tr{\A}<-2$ is not the exponential of any real matrix $X\in\sL(2)$. 
%
The logarithm on $\SL(2)$ can also be given in closed form \cite[(1.175)]{Tarantola06}.%
\footnote{In fact, the logarithm on diagonal matrices $D$ in $\SL(2)$ with positive eigenvalues is simple. Their trace is always $\lambda+\frac{1}{\lambda}\ge 2$ if $\lambda>0$. We infer that $\cosh(s)=\frac{\tr{D}}{2}$ can always be solved for $s$. Observe $\cosh(s)^2-\sinh(s)^2=1$.}
On the set $\SL(2)_\id$ the principal matrix logarithm is real and we have
\begin{align}
\label{logarithm_SL2}
\log &:\SL(2)_\id\mapsto\sL(2)_0\, ,\quad
   \log[S]:=\frac{s}{\sinh{s}}\left(S-\cosh(s)\id_2\right)\, ,\quad \cosh(s)=\frac{\tr{S}}{2} \, .
\end{align}

\subsubsection{Minimizing $\norm{\logp Q^T D}_F^2$ for $Q^T D\in \SL(2)_\id$}    
Let us define the open set
\begin{align}
  \mathcal{R}_D:=\{Q\in \SO(2)\, | \, Q^T D\in \SL(2)_\id \, \}\, .
\end{align}
On $\mathcal{R}_D$ the evaluation of $\log \overline{R}^T D$ is in the domain of the principal logarithm.
\noindent We are now able to show the optimality result with respect to rotations in $\mathcal{R}_D$. We use the given formula \eqref{logarithm_SL2} for the real logarithm on $\SL(2)_\id$ to compute
\begin{align}
\inf_{\overline{R}\in\mathcal{R}_D}& 
\norm{\log \overline{R}^TD }_F^2=
\inf_{\overline{R}\in\mathcal{R}_D} \frac{s^2}{\sinh(s)^2} \norm{\overline{R}^TD-\cosh(s)\id_2}_F^2\notag\\
&=\inf_{\overline{R}\in\mathcal{R}_D} \frac{s^2}{\sinh(s)^2} \left(\norm{\overline{R}^TD}_F^2-2\cosh(s)\Mprod{\overline{R}^TD}{\id_2}+2\cosh(s)^2\right)\notag\\
&\hspace{3cm} (\text{using}\;\cosh(s)=\frac{\tr{\overline{R}^TD}}{2})\notag\\
&=\inf_{\overline{R}\in\mathcal{R}_D} \frac{s^2}{\sinh(s)^2} \left(\norm{\overline{R}^TD}_F^2-4\cosh(s)^2+2\cosh(s)^2\right)\notag\\
&=\inf_{\overline{R}\in\mathcal{R}_D} \frac{s^2}{\sinh(s)^2} \left(\norm{\overline{R}^TD}_F^2-2\cosh(s)^2\right) \\
&=\inf_{\overline{R}\in\mathcal{R}_D} \frac{s^2}{\sinh(s)^2} \left(\norm{\overline{R}^TD}_F^2-\frac{1}{2}\Mprod{\overline{R}^TD}{\id}^2\right)\notag\\
&\ge \left( \inf_{\overline{R}\in\mathcal{R}_D\, , \cosh(s)=\frac{\tr{\overline{R}^TD}}{2}} \frac{s^2}{\sinh(s)^2} \right)\, 
\left(\inf_{\overline{R}\in\mathcal{R}_D} \left(\norm{\overline{R}^TD}_F^2-\frac{1}{2}\Mprod{\overline{R}^TD}{\id}^2\right)\right)\notag\\
&\hspace{6cm} \text{optimality of the orthogonal factor}\notag\\
&=\left(\inf_{\overline{R}\in\mathcal{R}_D\, , \cosh(s)=\frac{\tr{\overline{R}^TD}}{2}} \frac{s^2}{\sinh(s)^2} \right)\,
\left(\norm{D}_F^2-\frac{1}{2}\Mprod{D}{\id}^2\right)\notag\\
&=\norm{\dev_2 D}_F^2\, \inf_{\overline{R}\in\mathcal{R}_D\, , \cosh(s)=\frac{\tr{\overline{R}^TD}}{2}} \frac{s^2}{\sinh(s)^2}\notag
\end{align}
and on the other hand
\begin{align}
\label{log_computation}
  \norm{\log D}_F^2&=\frac{s^2}{\sinh(s)^2}\norm{D-\frac{1}{2}\tr{D}\id_2}_F^2\, ,\quad \quad\text{with $s\in\R$ s. that}\; \cosh(s)=\frac{\tr{D}}{2}=\frac{1}{2}\left(\lambda+\frac{1}{\lambda}\right)\notag\\
  &=\frac{s^2}{\sinh(s)^2}\norm{\dev_2 D}_F^2=\frac{s^2}{\sinh(s)^2}\, \frac{1}{2}(\lambda_1-\lambda_2)^2\quad\text{since}\;\norm{\dev_2 D}_F^2=\frac{1}{2}(\lambda_1-\lambda_2)^2\notag\\  
  &=\frac{s^2}{\sinh(s)^2}\, \frac{1}{2}(\lambda-\frac{1}\lambda)^2
  =\frac{s^2}{1+\sinh(s)^2-1}\, \frac{1}{2}(\lambda-\frac{1}\lambda)^2\notag\\   
   &=\frac{s^2}{\cosh(s)^2-1}\,\frac{1}{2}(\lambda-\frac{1}\lambda)^2  
  =\frac{s^2}{\frac{1}{4}(\lambda+\frac{1}{\lambda}   )^2-1}\, \frac{1}{2}(\lambda-\frac{1}\lambda)^2 \\
  &=\frac{(\arcosh(\frac{\lambda+\frac{1}{\lambda}}{2}))^2}{\frac{1}{4}(\lambda+\frac{1}{\lambda}   )^2-1}\, \frac{1}{2}(\lambda-\frac{1}\lambda)^2
  =\frac{(\arcosh(\frac{\lambda+\frac{1}{\lambda}}{2}))^2}{\frac{1}{4}(\lambda-\frac{1}{\lambda}   )^2}\,\frac{1}{2}(\lambda-\frac{1}\lambda)^2\notag\\  
  &=2 \, (\arcosh(\frac{\lambda+\frac{1}{\lambda}}{2}))^2=2 (\log\lambda)^2\, . \notag
  \end{align}
  The last equality can be seen by using the identity $\arcosh(x)=\log(x+\sqrt{x^2-1}),\, x\ge 1$ and setting $x=\frac{1}{2}\left(\lambda+\frac{1}{\lambda}\right)$, where we note that $x\ge 1$ for $\lambda>0$. Comparing the last result \eqref{log_computation} with the simple formula for the principal logarithm on $\SL(2)$ for diagonal matrices $D\in\SL(2)$ with positive eigenvalues shows
\begin{align}
    \norm{\log
    \begin{bmatrix}
        \lambda & 0\\
        0 & \frac{1}{\lambda}
        \end{bmatrix}}_F^2= \norm{
    \begin{bmatrix}
        \log\lambda & 0\\
        0 & \log \frac{1}{\lambda}
        \end{bmatrix}}_F^2
   = (\log\lambda)^2+(\log 1-\log\lambda)^2=2(\log\lambda)^2\, .
          \end{align}
To finalize the $\SL(2)$ case we need to show that
\begin{align}
&\inf_{\overline{R}\in\mathcal{R}_D\, , \cosh(s)=\frac{\tr{\overline{R}^TD}}{2}} \frac{s^2}{\sinh(s)^2}=
\frac{\hat{s}^2}{\sinh(\hat{s})^2}\, ,
\quad \text{with $\hat{s}\in\R$ s. that}\; \cosh(\hat{s})=\frac{\tr{D}}{2}=\frac{1}{2}\left(\lambda+\frac{1}{\lambda}\right)\, .
\end{align}
In order to do this, we write
\begin{align*}
\inf_{\overline{R}\in\mathcal{R}_D\, , \cosh(s)=\frac{\tr{\overline{R}^TD}}{2}} \frac{s^2}{\sinh(s)^2}
&=\inf_{\overline{R}\in\mathcal{R}_D\, , \cosh(s)=\frac{\tr{\overline{R}^TD}}{2}} \frac{s^2}{\cosh(s)^2-1}
=\inf_{\xi} \, \frac{\arcosh(\xi)^2}{\xi^2-1} \quad\quad \text{for}\; \xi=\frac{\tr{\overline{R}^TD}}{2}\, .
\end{align*}
One can check that the function 
\begin{align}
 g&: [1,\infty)\rightarrow \R^+\, ,\quad 
 g(\xi):=\frac{\arcosh(\xi)^2}{\xi^2-1}
\end{align}
is strictly monotone decreasing.
Thus 
$
g(\xi)=\frac{\arcosh(\xi)^2}{\xi^2-1}
$
is the smaller, the larger $\xi$ gets. The largest value for $\xi=\frac{\tr{\overline{R}^TD}}{2}$ is realized by $\hat{\xi}=\frac{\tr{D}}{2}$. Therefore
\begin{align*}
   \inf_{\xi}\,  \frac{\arcosh(\xi)^2}{\xi^2-1} \ge  \frac{\arcosh(\hat{\xi})^2}{\hat{\xi}^2-1}
   =\frac{\hat{s}^2}{\sinh(\hat{s})^2}\, ,\quad \quad\text{with $\hat{s}\in\R$ s. that}\; \cosh(\hat{s})=\frac{\tr{D}}{2}\, .\qeda
\end{align*}

\bibliographystyle{plain}
\bibliography{literatur1}


\def\noopsort#1{}\def\l{\char32l}\def\v#1{{\accent20 #1}}
  \let\^^_=\v\def\hbk{hardback}\def\pbk{paperback}

\end{document}